# Smooth classification of Cartan actions of higher rank semisimple Lie groups and their lattices

By Edward R. Goetze and Ralf J. Spatzier*

## Abstract

Let $G$ be a connected semisimple Lie group without compact factors whose real rank is at least 2, and let $\Gamma \subset G$ be an irreducible lattice. We provide a $C^\infty$ classification for volume-preserving Cartan actions of $\Gamma$ and $G$. Also, if $G$ has real rank at least 3, we provide a $C^\infty$ classification for volume-preserving, multiplicity free, trellised, Anosov actions on compact manifolds.

## 1. Introduction

Anosov diffeomorphisms and flows are some of the best understood and most important dynamical systems. They are the prototype of hyperbolic dynamical systems and enjoy special rigidity properties such as structural stability. Indeed, D. Anosov showed that a sufficiently small $C^1$ perturbation of an Anosov diffeomorphism is conjugate to the original diffeomorphism by a homeomorphism [1]. In this paper we will study Anosov actions of more general groups than $\mathbb{Z}$ and $\mathbb{R}$. By an Anosov action, we mean a locally faithful action of a (not necessarily connected) Lie group which contains an element which acts normally hyperbolically to the orbit foliation. This generalizes a definition of such actions by C. Pugh and M. Shub in [22]. Note that an Anosov action of a discrete group is simply an action such that some element of this group acts by an Anosov diffeomorphism.

Anosov actions of higher rank abelian or semisimple groups and their lattices are markedly different from Anosov diffeomorphisms and flows. In fact, during the last decade remarkable rigidity properties of actions of higher rank groups were discovered, ranging from local smooth rigidity to rigidity of invariant measures. Consider the standard action of $\mathrm{SL}(n, \mathbb{Z})$ on the $n$ torus, a prime example of an Anosov action of a lattice in a semisimple Lie group.

*The first author was supported in part by grants from the NSF and the University of Michigan. The second author was supported in part by a grant from the NSF.

In 1986, R. Zimmer conjectured that for $n > 2$, any sufficiently small $C^1$ perturbation of this action is *smoothly* conjugate to the standard action [30]. Infinitesimal, deformation and finally smooth local rigidity were established for this action in a sequence of papers by J. Lewis, S. Hurder, A. Katok and R. Zimmer [20], [10], [12], [16], [15], [17] and later generalized to other toral and nilmanifold actions by N. Qian [23], [24], [27].

Hurder actually conjectured that any Anosov action of a lattice in a higher rank semisimple group is essentially algebraic [10]. We will prove this conjecture for a special class of Anosov actions of lattices and a more general one for groups. The first are the *Cartan actions* introduced by Hurder in [10]. They are characterized by the property that suitable intersections of stable manifolds of certain commuting elements of the action are one dimensional (cf. Definition 3.8). The second class, also introduced by Hurder, is that of trellised actions. If $A \subset G$ is an abelian subgroup, then we call an Anosov action of $G$ *trellised with respect to $A$* if there exists a sufficiently large collection of one dimensional foliations invariant under the action of $A$ (cf. Definition 2.1). Cartan actions are always trellised. Finally, we will also use the notion of a multiplicity free action. These actions are characterized by the property that the super-rigidity homomorphism corresponding to the action consists of irreducible subrepresentations which are all multiplicity free (cf. Definition 3.3).

To clarify what we consider an essentially algebraic action we provide the following:

*Definition* 1.1. Let $H$ be a connected, simply connected Lie group with $\Lambda \subset H$ a cocompact lattice. Define $\mathrm{Aff}(H)$ to be the set of diffeomorphisms of $H$ which map right invariant vector fields on $H$ to right invariant vector fields. Define $\mathrm{Aff}(H/\Lambda)$ to be the diffeomorphisms of $H/\Lambda$ which lift to elements of $\mathrm{Aff}(H)$. Finally, we define an action $\rho : G \times H/\Lambda \to H/\Lambda$ to be *affine algebraic* if $\rho(g)$ is given by some homomorphism $\sigma : G \to \mathrm{Aff}(H/\Lambda)$.

THEOREM 1.2. *Let $G$ be a connected semisimple Lie group without compact factors and with real rank at least three, and let $A \subset G$ be a maximal $\mathbb{R}$-split Cartan subgroup. Let $M$ be a compact manifold without boundary, and let $\mu$ be a smooth volume form on $M$. If $\rho : G \times M \to M$ is an Anosov action on $M$ which preserves $\mu$, is multiplicity free, and is trellised with respect to $A$, then, by possibly passing to a finite cover of $M$, $\rho$ is $C^\infty$ conjugate to an affine algebraic action, i.e., there exist*

1. *a finite cover $M' \to M$,*

2. *a connected, simply connected Lie group $L$,*

3. *a cocompact lattice $\Lambda \subset L$,*

4. *a $C^\infty$ diffeomorphism $\phi : M \to L/\Lambda$, and*

5. *a homomorphism $\sigma : G \to \operatorname{Aff}(L/\Lambda)$*

*such that $\rho'(g) = \phi^{-1}\sigma(g)\phi$, where $\rho'$ denotes the lift of $\rho$ to $M'$.*

If, for a given Cartan subgroup $A \subset G$, the nontrivial elements of the Oseledec decomposition of $TM = \oplus E_i$ corresponding to $A$ consist entirely of one dimensional spaces, then it follows that the action must be both trellised and multiplicity free. This yields the following:

COROLLARY 1.3. *Let $G$ be a connected semisimple Lie group without compact factors and with real rank at least three, and let $A \subset G$ be a maximal $\mathbb{R}$-split Cartan subgroup. Let $M$ be a compact manifold without boundary, and let $\mu$ be a smooth volume form on $M$. If $\rho : G \times M \to M$ is an Anosov action on $M$ which preserves $\mu$ and is such that the nontrivial elements of the Oseledec decomposition with respect to $A$ consists of one dimensional Lyapunov spaces, then, by possibly passing to a finite cover of $M$, $\rho$ is $C^\infty$ conjugate to an affine algebraic action.*

We obtain $\mathbb{R}$-rank 2 results with an additional assumption.

COROLLARY 1.4. *Assume the conditions of Theorem 1.2. If, in addition, the trellis consists of one dimensional strongest stable foliations, i.e. the action is Cartan, then the above classification holds when the real rank of $G$ is at least two.*

The next results provide a similar classification for actions of lattices.

THEOREM 1.5. *Let $G$ be a connected semisimple Lie group without compact factors such that each simple factor has real rank at least 2, and let $\Gamma \subset G$ be a lattice. Let $M$ be a compact manifold without boundary and $\mu$ a smooth volume form on $M$. Let $\rho : \Gamma \times M \to M$ be a volume-preserving Cartan action. Then, on a subgroup of finite index, $\rho$ is $C^\infty$ conjugate to an affine algebraic action.*

*More specifically, on a subgroup of finite index, $\rho$ lifts to an action of a finite cover $M' \to M$ which is $C^\infty$ conjugate to the standard algebraic action on the nilmanifold $\tilde{\pi}_1(M')/\pi_1(M')$, where $\tilde{\pi}_1(M')$ denotes the Malcev completion of the fundamental group of $M'$, i.e., the unique, simply connected, nilpotent Lie group containing $\pi_1(M')$ as a cocompact lattice.*

We point out that Theorem 1.5 proves Hurder's conjecture of Anosov rigidity of lattice actions in the case of Cartan actions [10].

As an immediate corollary, we recover the local rigidity results obtained for Cartan homogeneous actions.

COROLLARY 1.6. *Let $\Gamma \subset G$ be an irreducible lattice as in Theorem* 1.5, *and let $\phi : \Gamma \times M \to M$ be a volume-preserving Cartan action on a closed manifold $M$. Then $\phi$ is locally $C^\infty$ rigid.*

Rigidity of higher rank groups and their actions is typically connected with an analysis of the action of a maximal abelian subgroup $A$ of the original group. As a first step in the proof we show that there always exists a Hölder Riemannian metric on the manifold with respect to which $A$ has uniform expansion and contraction. For $G$ actions, we proved this in an earlier paper [7]. For lattices, this follows from a result of N. Qian on the existence of a continuous framing which transforms under $G$ according to some finite dimensional representation of $G$ [25].

The main contribution of the current paper is an analysis of the regularity of this metric and of various unions of stable and unstable foliations. This analysis involves only the abelian subgroup $A$. In fact, in Section 2 we present an abstract version of this for general trellised Anosov actions of $\mathbb{R}^k$. A key ingredient of the argument is the construction of isometries of subfoliations of the manifold using an element of $A$ which does not expand or contract the leaves. This is an idea due to A. Katok and was employed in [18] to control invariant measures for hyperbolic actions of higher rank abelian groups.

In Section 3, we consider the semisimple situation. At this point, we have a smooth framing of the manifold which transforms according to a finite dimensional representation of $G$. We then adapt an argument of G. A. Margulis and N. Qian [25] to finish the proof of our main results.

We thank G. Prasad, C. Pugh, F. Raymond, and M. Brown for several helpful discussions.

## 2. Smooth geometric structures for $\mathbb{R}^k$ actions

In this section, we consider a certain class of $\mathbb{R}^k$ actions on a closed manifold $M$ with constant derivative with respect to some Hölder framing. By analyzing the behavior of this action on certain stable and unstable subfoliations, we show that this framing is actually smooth. This is the key ingredient in the classification of the actions considered in Section 3.

2.1. *Preliminaries.* We shall assume that $A = \mathbb{R}^k$ acts smoothly on a closed manifold $M$ preserving a measure $\mu$. For any $a \in A$, we have a Lyapunov decomposition of the tangent bundle, with Lyapunov exponents $\{\chi_i\}$. Since $A$ is abelian, we may refine this decomposition to a joint splitting $TM = \bigoplus E_i$ for all $a \in A$. Note that the exponents still vary with the choice of $a \in A$. Because $A$ is abelian, and may be identified with its Lie algebra, we can consider the exponents as linear functionals on $A$, which, henceforth, will be referred to as the *weights* of the action with respect to $\mu$. Let $\mathcal{W}(A)$ denote the set of such weights for this action.

We present a modified version of Hurder's definition of a trellised action. We will call two foliations *pairwise transverse* if their tangent spaces intersect trivially. The standard notion in differential topology also requires the sum of the tangent spaces to span the tangent space of the manifold. This condition is replaced by the first condition in the definition below.

*Definition* 2.1. Let $A$ be an (abelian) group. A $C^\infty$ action $\phi: A \times X \to X$ is *trellised* if there exists a collection $\mathcal{T}$ of one dimensional, pairwise transverse foliations $\{\mathcal{F}_i\}$ of $X$ such that

1. The tangential distributions have internal direct sum
$$T\mathcal{F}_1 \oplus \cdots \oplus T\mathcal{F}_r \oplus TA \cong TX,$$
where $TA$ is the distribution tangent to the $A$ orbit.

2. For each $x \in X$ the leaf $L_i(x)$ of $\mathcal{F}_i$ through $x$ is a $C^\infty$ immersed submanifold of $X$.

3. The $C^\infty$ immersions $L_i(x) \to X$ depend uniformly Hölder continuously on the basepoint $x$ in the $C^\infty$ topology on immersions.

4. Each $\mathcal{F}_i$ is invariant under $\phi(a)$ for every $a \in A$.

Moreover, if a group $H$ acts on a manifold $M$ and $A \subset H$ is an abelian subgroup, then we say the action is *trellised with respect to $A$* if the action restricted to $A$ is trellised.

Later in this paper, we will consider the case where $H$ is a semisimple Lie group without compact factors and $A$ is a maximal $\mathbb{R}$-split Cartan subgroup.

*Example* 2.2 (trellised actions). 1. Let $G = \mathrm{SO}(n,n)$, the $\mathbb{R}$-split group with Lie algebra $\mathfrak{b}_n$, and let $M = \mathrm{SO}(n, n+1)/\Lambda$ for some cocompact lattice $\Lambda \subset \mathrm{SO}(n, n+1)$. Suppose that the action of $G$ on $M$ comes from the standard inclusion $\mathrm{SO}(n,n) \hookrightarrow \mathrm{SO}(n, n+1)$. The set of weights for this inclusion is the union of the roots for $\mathfrak{b}_n$ and the weights corresponding to the standard action of $\mathrm{SO}(n, n+1)$ on $\mathbb{R}^{2n+1}$. In particular, each weight space is one dimensional and no weight is a positive multiple of any other. It follows that this action is trellised, and all the nontrivial Lyapunov spaces are one dimensional. However, it is not Cartan (cf. Definition 3.8), since the weight spaces corresponding to weights of the standard action cannot be written as the strongest stable space for any element in $\mathrm{SO}(n,n)$.

2. For simpler (transitive) examples, consider an $\mathbb{R}$-split semisimple connected Lie group $G$ without compact factors. If $\Lambda \subset G$ is a cocompact lattice, then the natural $G$ action on $G/\Lambda$ will be trellised.

Let us return to the case of an action of an abelian group $A$ on a compact manifold $M$. Throughout this section, we shall make the following assumptions:

(A0) The action is locally free.

(A1) The Lyapunov decomposition extends to a Hölder splitting $TM = \bigoplus E_i$ of the tangent bundle.

(A2) There exists an $A$-invariant smooth volume on $M$.

(A3) The action of any 1 parameter subgroup of $A$ is ergodic on $M$.

(A4) The action on $M$ is trellised with respect to $A$.

(A5) There exists a Hölder Riemannian metric on $M$ such that $\|av\| = e^{\chi_i(a)}\|v\|$ for every $a \in A$ and for every $v \in E_i$.

(A6) If $E_i \not\subset TA$, then $\chi_i \not\equiv 0$.

Since for an ergodic flow $\{\phi_t\}_{t\in\mathbb{R}}$, the map $\phi_{t_0}$ is ergodic for almost every $t_0 \in \mathbb{R}$, we can replace Assumption (A3) with the equivalent assumption:

(A3′) Every 1 parameter subgroup of $A$ contains an ergodic element.

An immediate consequence of these assumptions is that the $A$ action on $M$ is Anosov, i.e., there exists some element in $A$ that acts normally hyperbolically on $M$ with respect to the $A$ foliation. In fact, every element in the complement of the union of the hyperplanes $\ker(\chi_i)$, for $i$ such that $E_i \not\subset TA$, is normally hyperbolic. We also point out another immediate consequence.

LEMMA 2.3. *The trellis is subordinate to the Lyapunov decomposition, i.e., for every $i$, there exists some $j$ such that $T\mathcal{F}_i \subset E_j$.*

*Proof.* Suppose that $T\mathcal{F}_i \subset \bigoplus_{j \in \mathcal{J}} E_j$ where $\mathcal{J}$ is the smallest possible set of indices. Pick an ergodic element $a \in A$. As $n \to \infty$, $da^n(T\mathcal{F}_i)$ converges into $E_{j_1}$ where $\chi_{j_1}(a)$ is the maximum value of $\{\chi_j(a)\}_\mathcal{J}$. Similarly, as $n \to -\infty$, $da^n(T\mathcal{F}_i)$ converges into $E_{j_2}$ where $\chi_{j_2}(a)$ is the minimum value of $\{\chi_j(a)\}_\mathcal{J}$. Because of the recurrence property of the action, continuity of the $\mathcal{F}_i$, and the assumption that $\mathcal{F}_i$ is fixed by $A$, we are presented with a contradiction unless $j_1 = j_2$, i.e., unless $T\mathcal{F}_i \subset E_j$ for some $j$. □

The main result of this section is that the geometric structures on $M$ have significantly greater regularity than initially assumed.

THEOREM 2.4. *Suppose $A = \mathbb{R}^k, k \geq 3$, acts on a closed manifold $M$ satisfying Assumptions* (A0) *through* (A6). *Then the trellis $\mathcal{T}$ and the Riemannian metric in* (A5) *are both $C^\infty$. In particular, the $C^\infty$ immersions $L_i(x) \to M$ depend $C^\infty$ on the basepoint $x$ in the $C^\infty$ topology on immersions, and each $\mathcal{F}_i$ has uniformly $C^\infty$ leaves.*

Since the proof proceeds through a number of steps, we provide a brief outline. First, we define a distribution $\mathcal{N}_\mathcal{H}^+$ of $TM$ consisting of a particular collection of stable directions and show that it is an integrable distribution tangent to a Hölder foliation with $C^\infty$ leaves $N_\mathcal{H}^+(x)$. By restricting the Hölder metric on $M$ to the leaves of this foliation, we can consider the group of isometries of a particular leaf. We then show that there exists a subgroup of isometries that acts simply transitively on $N_\mathcal{H}^+(x)$. The idea is that certain elements $a \in A$ as well as limits of certain sequences of powers of such an $a$ are isometries between the leaves of $N_\mathcal{H}^+$. We continue by showing that there exists a canonically defined set of these limiting isometries which acts simply transitively on $N_\mathcal{H}^+(x)$.

The second step is to consider a larger foliation $N_\mathcal{H}$ of $M$ with leaves that consist both of certain stable and unstable directions. We define a new metric on $N_\mathcal{H}(x)$, and show that its group of isometries acts transitively. Using Montgomery and Zippin's work on Hilbert's Fifth Problem, we conclude that $N_\mathcal{H}(x)$ is a homogeneous space of a Lie group. This yields a new differentiable structure on $N_\mathcal{H}(x)$ for all $x$ with respect to which the part of the trellis tangent to $N_\mathcal{H}(x)$ is automatically smooth on $N_\mathcal{H}(x)$.

The final step in the proof of Theorem 2.4 is to show that $N_\mathcal{H}(x)$ with its differentiable structure as a homogeneous space smoothly immerses into $M$ via the inclusion $N_\mathcal{H}(x) \hookrightarrow M$. Theorem 2.4 then quickly follows. To accomplish this, we use an argument similar to that presented by Katok and Lewis in [16] where they use the nonstationary Sternberg linearization to show that an *a priori* topological conjugacy is actually smooth. We note that in light of our assumptions, we require only a simplified version of Katok and Lewis' original argument.

2.2. *Simply transitive groups of isometries for stable subfoliations.* Fix some $b_0 \in A$ once and for all. Suppose $\mathcal{H} \subset A$ is a proper vector subspace. Define $\mathcal{J}_\mathcal{H}^+ = \{\chi_i \in \mathcal{W}(A) | \chi_i(b_0) > 0 \text{ and } \mathcal{H} \subset \ker(\chi_i)\}$, and set

$$\mathcal{N}_\mathcal{H}^+ = \bigoplus_{\chi_i \in \mathcal{J}_\mathcal{H}^+} E_i.$$

We can similarly define $\mathcal{N}_\mathcal{H}^-$. Of course, most interesting is the case where $\mathcal{J}_\mathcal{H}^+$ is not empty. We shall call kernels of nonzero weights $\chi_i$ *weight hyperplanes*. Since there are only finitely many weights, there are also only finitely many weight hyperplanes.

If $a \in A$ does not lie on any of the weight hyperplanes, then $a$ is a *normally hyperbolic* or *regular element*. If $\mathcal{H}$ is a weight hyperplane then call $a \in \mathcal{H}$ *generic*, if for every weight $\chi$, $\chi(a) = 0$ implies $\chi(\mathcal{H}) = 0$. For $a \in A$, let

$$E_a^+ = \bigoplus_{\{\chi_i \in \mathcal{W}(A), \chi_i(a) > 0\}} E_i.$$

We similarly define $E_a^-$ and $E_a^0$.

LEMMA 2.5.
 Let $\mathcal{P} \subset \mathcal{W}(A)$ be the set of weights which are positive on $b_0$. Then there exist

1. *an ordering of* $\mathcal{P} = \{\chi_1, \ldots, \chi_r\}$, *and*

2. *regular elements* $b_i \in A, 1 \leq i \leq r$,

such that $E_{b_0}^+ \cap E_{b_i}^+ = \bigoplus_{j=1}^i E_j$. Hence, $\bigoplus_{j=1}^i E_j$ *forms an integrable distribution tangent to a Hölder foliation with uniformly $C^\infty$ leaves.*

*Proof.* Let $P$ be a two dimensional plane in $A$ which contains $b_0$. Then $P$ is not contained in any weight hyperplane. Thus, the intersection of any weight hyperplane with $P$ is a one dimensional line. Let $\mathcal{L}_\chi = P \cap \ker(\chi)$ for every $\chi \in \mathcal{P}$. These lines divide $P$ into $2r$ distinct sectors such that $\pm b_0 \notin \mathcal{L}_\chi$ for every $\chi \in \mathcal{P}$. Let $-B_0$ be the region in $P$ which contains $-b_0$, and pick $B_1$ to be a region adjacent to $-B_0$. Pick $b_1 \in B_1$ to be some regular element. For every $1 < i \leq r$, let $B_i$ be the unique region adjacent to $B_{i-1}$ not equal to $B_{i-2}$ (or $-B_0$ if $i = 2$), and pick $b_i \in B_i$ to be some regular element. Note that $B_r$ contains $b_0$ so that we may choose $b_r = b_0$.

Let $\chi_i$ be the element of $\mathcal{P}$ such that $P \cap \ker(\chi_i)$ separates $B_{i-1}$ and $B_i$ ($-B_0$ and $B_1$ when $i = 1$). It follows that $b_0$ and $b_i$ lie on the same side of $\mathcal{L}_{\chi_j}$ whenever $j \leq i$, and on opposite sides of $\mathcal{L}_{\chi_j}$ whenever $j > i$. We therefore conclude $E_{b_0}^+ \cap E_{b_i}^+ = \bigoplus_{j=1}^i E_j$. The final comment follows from [28, App. IV, Th. IV.1]. □

*Remark* 2.6. This proof easily generalizes to produce an ordering of the weights in $\mathcal{J}_\mathcal{H}^+$, and regular elements $\{b_i\} \in A$ such that $E_{b_0}^+ \cap E_{b_i}^+ \cap \mathcal{N}_\mathcal{H}^+ = \bigoplus_{j=1}^i E_j$. In conjunction with the following lemma, we have that $\left(\bigoplus_{j=1}^i E_{b_j}^+\right) \cap \mathcal{N}_\mathcal{H}^+$ is an integrable distribution tangent to a Hölder foliation with uniformly $C^\infty$ leaves. Hence, we can produce a nested sequence of Hölder foliations $\mathcal{L}_1 \subset \mathcal{L}_2 \subset \cdots \subset \mathcal{N}_\mathcal{H}^+$ with $C^\infty$ leaves such that $\left(\bigoplus_{j=1}^i E_{b_j}^+\right) \cap \mathcal{N}_\mathcal{H}^+$ is the distribution tangent to $\mathcal{L}_i$.

LEMMA 2.7. *Suppose $\mathcal{H} \subset A$ is a proper linear subspace contained in some weight hyperplane. There exist regular elements $c, d \in A$ such that $\mathcal{N}_\mathcal{H}^+ = E_c^+ \cap E_d^-$. Hence, $\mathcal{N}_\mathcal{H}^+(x)$ is an integrable distribution tangent to the intersection $N_\mathcal{H}^+(x) = W_c^+(x) \cap W_d^-(x)$, which forms a Hölder foliation with $C^\infty$ leaves. A similar result holds for $\mathcal{N}_\mathcal{H}^-$.*

*Proof.* Let $P$ be a two dimensional plane containing $b_0$ and some nonzero $a \in \mathcal{H}$. Let $c = a - \varepsilon b_0$ and $d = a + \varepsilon b_0$. If $\varepsilon$ is small enough then the only weight hyperplane that the line segment from $c$ to $d$ intersects contains $\mathcal{H}$. In particular, we may assume that $c$ and $d$ are regular.

If $\mu \in \mathcal{W}(A)$ and $\mu(c)$ and $\mu(d)$ are both greater than 0, we must have $\mu(a) > 0$ since $a = (c+d)/2$. By choice of $\varepsilon$, if $\mu(a) > 0$, then $\mu(c)$ and $\mu(d)$ are both greater than 0. It follows that $E_a^+ = E_c^+ \cap E_d^+$, and similarly that $E_a^- = E_c^- \cap E_d^-$. In fact, an analogous argument shows that we can write $E_a^0 = \mathcal{N}_\mathcal{H}^+ + \mathcal{N}_\mathcal{H}^- + T(A)$ or $E_a^0 = (E_c^+ \cap E_d^-) + (E_c^- \cap E_d^+) + T(A)$ where $T(A)$ represents the tangent space to the orbit. It follows that $\mathcal{N}_\mathcal{H}^+ = E_c^- \cap E_d^+ = E_{-c}^+ \cap E_d^+$. $\square$

Consider the restriction of the Hölder metric on $M$ to $N_\mathcal{H}^+(x)$, and let $\mathrm{Isom}(N_\mathcal{H}^+(x))$ be the set of isometries with respect to this metric. By [3], any isometry with respect to this metric must be at least $C^1$. Denote by $\hat{\mathcal{I}}(N_\mathcal{H}^+(x))$ the subgroup of $\mathrm{Isom}(N_\mathcal{H}^+(x))$ which preserves the tangent bundle for each element of the trellis belonging to $\mathcal{N}_\mathcal{H}^+(x)$, i.e., if $\phi \in \hat{\mathcal{I}}(N_\mathcal{H}^+(x))$, then $d\phi_x(T\mathcal{F}_i(x)) = T\mathcal{F}_i(\phi(x))$ for every $\mathcal{F}_i(x) \subset N_\mathcal{H}^+(x)$. Let $\mathcal{I}(N_\mathcal{H}^+(x))$ be the connected component of the identity of $\hat{\mathcal{I}}(N_\mathcal{H}^+(x))$.

THEOREM 2.8. *Let $\mathcal{H} \subset A$ be a proper linear subspace contained in some weight hyperplane. Then $\mathcal{I}(N_\mathcal{H}^+(x))$ acts simply transitively on $N_\mathcal{H}^+(x)$ for every $x \in M$.*

The main step in the proof is to demonstrate the existence of a certain class of isometries.

PROPOSITION 2.9. *Let $a \in \mathcal{H} \subset A$ and suppose $\{n_k\}$ is a sequence such that $\lim_{k \to \infty} a^{n_k} x = y$. Then there exist a subsequence $\{m_j\}$ and a map $\alpha : N_\mathcal{H}^+(x) \to N_\mathcal{H}^+(y)$ such that $\alpha(z) = \lim_{j \to \infty} a^{m_j}(z)$ and $\alpha$ is an isometry with respect to the relevant induced Hölder metrics.*

We need a few basic lemmas.

LEMMA 2.10. *Let $\langle \cdot, \cdot \rangle_\infty$ be a $C^\infty$ Riemannian metric on $M$ and consider its restriction to $N_\mathcal{H}^+(x)$. Let $\exp_x : \mathcal{N}_\mathcal{H}^+(x) \to N_\mathcal{H}^+(x)$ denote the corresponding exponential map. Then the map $x \mapsto \exp_x$ is $C^0$ in the $C^k$ topology; i.e., if*

$\phi : \mathbb{R}^l \times T \to M$ *is a local trivialization for the foliation, then the composition*

$$\mathbb{R}^l \stackrel{d\phi}{\to} \mathcal{N}_\mathcal{H}^+(x) \stackrel{\exp_x}{\to} N_\mathcal{H}^+(x) \stackrel{\phi^{-1}}{\to} \mathbb{R}^l \times \{x\} \stackrel{\text{proj}}{\to} \mathbb{R}^l$$

*depends $C^0$ in the $C^k$ topology on $x$.*

*Proof.* Note that $N_H^+(q)$ varies $C^0$ in the $C^k$ topology since stable manifolds vary $C^0$ in the $C^k$ topology and $N_H^+(x)$ is a transverse intersection of stable manifolds. Also note that $g_{ij} = \left\langle \frac{\partial}{\partial x_i}, \frac{\partial}{\partial x_j} \right\rangle_\infty$ is $C^0$ in $q$ which is a $C^{k-1}$ function on each $N_H^+(q)$. Choose an embedding $q : D^m \to M$, and pull back the metric on $M$ to a metric on $D^m$: $(D^m, g_q)$. Then $g_q$ is a $C^\infty$ metric on $D^m$ which varies $C^0$ in $q$ in the $C^k$ topology. This implies that the Christoffel symbols $\Gamma_{ij}^k$ vary continuously in $q$. The exponential map is the solution to a differential equation whose parameters vary continuously in $q$ since the $\Gamma_{ij}^k$ do. This implies that the solutions vary $C^0$ in $q$. Hence the exponential maps vary $C^0$ in $q$. □

The next lemma is an immediate corollary.

LEMMA 2.11. *Let $\langle \cdot, \cdot \rangle_\infty$ be a $C^\infty$ Riemannian metric on $M$. There exists a lower bound $\iota$ for the injectivity radius of $\langle \cdot, \cdot \rangle_\infty |_{N_\mathcal{H}^+(x)}$ which is independent of $x$.*

*Proof.* We will need the following slight generalization of the implicit function theorem, really a parametrized version of the inverse function theorem. We indicate a proof as we were unable to find a reference.

PROPOSITION 2.12. *Let $U$ be open in $\mathbb{R}^n$, $V$ open in $\mathbb{R}^m$, and let $F : U \times V \to \mathbb{R}^n$ be a map such that every restriction $F_v := F |_{U \times \{v\}}$ is $C^1$ on $U$ with derivative $f'_v = \text{Id}$. Assume further that the map $v \mapsto F_v$ is continuous in the $C^1$-topology. Then there exist open sets $U' \subset U$ and $V' \subset V$ such that for all $v \in V$, $F_v$ is a diffeomorphism from $U' \times \{v\}$ onto its image.*

*Proof.* Since the $F_v$ depend continuously in the $C^1$-topology, this follows straight away from the following standard estimate (cf. [19, p. 124]) of the size of the radius of a ball on which the maps $F_v$ are diffeomorphisms:

Consider a closed ball $\bar{B}_r(0) \subset U$ and a number $0 < s < 1$ such that $|F'_v(z) - F'_v(x)| \leq s$ for all $x, z \in B_r(0)$. Then $F_v$ is a diffeomorphism of $B_{r(1-s)}(0)$ onto its image. □

Applying this proposition to the exponential maps $\exp_x : \mathcal{N}_\mathcal{H}^+(x) \to N_\mathcal{H}^+(x)$, we see that for all $x$ there is a neighborhood $U_x$ of $x$ such that the $\exp_x$ are diffeomorphisms on balls about 0 in $\mathcal{N}_\mathcal{H}^+(x)$ of fixed radius. Covering the compact manifold $M$ with finitely many such neighborhoods $U_x$ shows that the injectivity radius is bounded below. □

*Proof of Proposition* 2.9. Let $\langle \cdot, \cdot \rangle$ and $d_x$ be the induced Hölder metric and corresponding distance function on $N_{\mathcal{H}}^+(x)$. Assume for the time being that for $i = 1, 2$, there exist $x_i \in N_{\mathcal{H}}^+(x), y_i \in N_{\mathcal{H}}^+(y)$ such that $\lim_{k \to \infty} a^{n_k} x_i = y_i$. Pick a $C^\infty$ Riemannian metric $\langle \cdot, \cdot \rangle_\infty$ on $M$ such that there exist constants $s$ and $S$ such that

$$s < \sqrt{\frac{\langle v, v \rangle}{\langle v, v \rangle_\infty}} < S$$

for any $v \in TM$. Thus, if $c$ is any curve in $N_{\mathcal{H}}^+(x)$ between $x_1$ and $x_2$,

$$s \cdot l_\infty(c) < l(c) < S \cdot l_\infty(c),$$

where $l$ and $l_\infty$ are the Hölder and $C^\infty$ lengths for the curve $c$ respectively.

Let $\iota$ be the bound on the injectivity radius obtained in Lemma 2.11. Pick $\varepsilon > 0$ and suppose $d(x_1, x_2) < \frac{s \cdot \iota}{1+\varepsilon}$. Let $c$ be a curve in $N_{\mathcal{H}}^+(x)$ from $x_1$ to $x_2$ such that $l(c) < (1+\varepsilon)d(x_1, x_2)$. Since $a \in \mathcal{H}$ implies $\chi_i(a) = 0$ for every $\chi_i \in \mathcal{J}_{\mathcal{H}}^+$, by Assumption (A5), we must have that $l(a(c)) = l(c)$, and hence $l(a^{n_k}(c)) = l(c) < (1+\varepsilon)d(x_1, x_2)$. Thus

$$l_\infty(a^{n_k}(c)) \leq \frac{1}{s} l(a^{n_k}(c)) < \frac{1+\varepsilon}{s} d(x_1, x_2) < \iota.$$

This implies that for every $k$ there exists a $v_k \in \mathcal{N}_{\mathcal{H}}^+(a^{n_k}(x_1))$ such that $\exp_{a^{n_k}(x_1)}(v_k) = a^{n_k}(x_2)$ and

$$\|v_k\|_\infty = d_\infty(a^{n_k}(x_1), a^{n_k}(x_2)) < \frac{1+\varepsilon}{s} d(x_1, x_2).$$

Pick a subsequence such that $v_k \to v \in \mathcal{N}_{\mathcal{H}}^+(y_1)$. By Lemma 2.10,

$$\exp_{y_1}(v) = \lim_{k \to \infty} \exp_{a^{n_k} x_1}(v_k) = \lim_{k \to \infty} a^{n_k} x_2 = y_2.$$

Thus,

$$d_\infty(y_1, y_2) \leq \|v\|_\infty \leq \frac{1+\varepsilon}{s} d(x_1, x_2).$$

But

$$\|v\|_\infty \leq \frac{1}{S} l(\exp_{y_1}(tv|_{t \in [0,1]})) \leq \frac{1}{S} d(y_1, y_2).$$

Consequently, for any $\varepsilon > 0$ and any $x_1, x_2 \in N_{\mathcal{H}}^+(x)$ such that $d(x_1, x_2) < \frac{s \cdot \iota}{1+\varepsilon}$,

(1) $$d(y_1, y_2) \leq \frac{S}{s} d(x_1, x_2).$$

However, for arbitrary $x_1, x_2 \in N_{\mathcal{H}}^+(x)$, we can divide any curve between $x_1$ and $x_2$ into a finite number of pieces each with length less than $\frac{s \cdot \iota}{1+\varepsilon}$. As a result, Equation 1 holds for any $x_1, x_2 \in N_{\mathcal{H}}^+(x)$. By choosing a $C^\infty$ Riemannian metric which better approximates the Hölder metric, we can ensure $\frac{S}{s}$ is arbitrarily close to 1. Thus if $\lim_{k \to \infty} a^{n_k}(x_i) = y_i$ for $i = 1, 2$, then

(2) $$d(y_1, y_2) \leq d(x_1, x_2).$$

Choose $\{x_i\}$ to be a countable dense subset of $N_{\mathcal{H}}^+(x)$. Since any $a \in A$ preserves the Lyapunov decomposition, it follows that $a$ maps $N_{\mathcal{H}}^+(x)$ to $N_{\mathcal{H}}^+(ax)$. Hence, setting $\{n_k^0\} = \{n_k\}$ then, using Equation 2 and compactness of $M$, for every $i$ there exists a subsequence $\{n_l^i\}$ of $\{n_k^{i-1}\}$ such that $a^{n_l^i}(x_i) \to y_i$ for some $y_i \in N_{\mathcal{H}}^+(y)$. Using a standard diagonal argument, there exists a subsequence $\{m_j\}$ of $\{n_k\}$ such that $a^{m_j}(x_i) \to y_i$ for every $i$. Define $\alpha(x_i) = y_i$. By Equation 2, we can extend $\alpha$ continuously to be defined on all of $N_{\mathcal{H}}^+(x)$ and such that $\alpha(z) = \lim_{k \to \infty} a^{n_k} z$ for all $z \in N_{\mathcal{H}}^+(x)$.

Summarizing, we have $\alpha : N_{\mathcal{H}}^+(x) \to N_{\mathcal{H}}^+(y)$ which is Lipschitz with Lipschitz constant $\leq 1$. To complete the proof, we will show that $\alpha$ has an inverse which is also Lipschitz with Lipschitz constant $\leq 1$.

Suppose $x_1, x_2 \in N_{\mathcal{H}}^+(x)$ and $d_\infty(x_1, x_2) < \iota$; i.e., $x_1$ and $x_2$ are within the bounds for the injectivity radius of the $C^\infty$ metric. Let $v \in T_{\alpha(x_1)} \mathcal{N}_{\mathcal{H}}^+(y)$ such that $\exp_{\alpha(x_1)}(v) = \alpha(x_2)$. Pick a sequence of vectors $v_k \in T_{a^{n_k}(x_1)} \mathcal{N}_{\mathcal{H}}^+(a^{n_k}(x_1))$ such that $\exp_{a^{n_k}(x_1)}(v_k) = a^{n_k}(x_2)$. By Lemma 2.10, $\exp_{a^{n_k}(x_1)}(v_k) \to \alpha(x_2)$, and by uniqueness of geodesics below the injectivity radius, we have that the limiting curve must be the geodesic $\exp_{\alpha(x_1)}(tv)$. Hence, we get $\exp_{a^{n_k}(x_1)}(tv_k) \to \exp_{\alpha(x_1)}(tv)$. As a result, $v_k \to v$, and we may conclude

$$d_\infty(\alpha(x_1), \alpha(x_2)) = \lim_{k \to \infty} d_\infty(a^{n_k}(x_1), a^{n_k}(x_2)).$$

Therefore, for any $x_1, x_2 \in N_{\mathcal{H}}^+(x)$ such that $d_\infty(x_1, x_2) < \iota$,

$$d(\alpha(x_1), \alpha(x_2)) \geq s \cdot d_\infty(\alpha(x_1), \alpha(x_2)) = s \cdot \lim_{k \to \infty} d_\infty(a^{n_k}(x_1), a^{n_k}(x_2))$$

$$\geq \frac{s}{S} \lim_{k \to \infty} d(a^{n_k}(x_1), a^{n_k}(x_2)) = \frac{s}{S} d(x_1, x_2).$$

When $\varepsilon = \frac{\iota}{2S}$, then if $d(x_1, x_2) < \varepsilon$, we have $d_\infty(x_1, x_2) < \varepsilon \cdot S \leq \frac{\iota}{2}$, and therefore $d(\alpha x_1, \alpha x_2) \geq \frac{s}{S} d(x_1, x_2)$. Now fix $\delta > 0$, and pick a $C^\infty$ Riemannian metric so that $\frac{s}{S} > 1 - \delta$. The argument above shows that there exists an $\varepsilon > 0$ such that

(3) $\quad d(x_1, x_2) < \varepsilon$ implies $d(\alpha(x_1), \alpha(x_2)) \geq (1 - \delta) d(x_1, x_2)$.

In particular, this shows that $\alpha$ is locally injective: If $d(x_1, x_2) < \varepsilon$, then $\alpha x_1 \neq \alpha x_2$. By Invariance of Domain [8, Cor. 18-9], $\alpha$ is an open map and is therefore a local homeomorphism.

Let $B_r(x)$ and $S_r(x)$ be the $r$-ball and the $r$-sphere in the Hölder metric about $x$, and let $\zeta = \varepsilon(1 - \delta)/3$. We claim

(4) $\quad\quad\quad\quad\quad\quad B_\zeta(\alpha(x_1)) \subset \alpha(B_{\varepsilon/2}(x_1))$

for any $x_1 \in N_{\mathcal{H}}^+(x)$. Suppose $y \in B_\zeta(\alpha(x_1))$, and let $\gamma(t), t \in [0, 1]$, be a path from $\alpha(x_1)$ to $y$ lying inside $B_\zeta(\alpha(x_1))$. Since $\alpha(x_1) = \gamma(0)$, the set $\{t | \gamma(t) \in \alpha(B_{\varepsilon/2}(x_1))\}$ is nonempty. Let $t_0$ be the supremum of this set. Since $\alpha$ is a local homeomorphism, $t_0 > 0$. Pick $t_n \to t_0$ and $x_n \in B_{\varepsilon/2}(x_1)$ such that $\gamma(t_n) = \alpha(x_n)$. Passing to a subsequence, we may assume that $x_n \to x'$ as $n \to \infty$. Then $\gamma(t_0) = \alpha(x')$. By Equation 3, we have $B_\zeta(\alpha(x_1)) \cap \alpha(S_{\varepsilon/2}(x_1))$ is empty. Hence $x' \in B_{\varepsilon/2}(x_1)$. Since $\alpha$ is a local homeomorphism, this yields a contradiction unless $t_0 = 1$. This proves Equation 4.

Since Equation 4 holds for all $x_1 \in N_{\mathcal{H}}^+(x)$, it follows $\alpha$ is a closed map. Since $N_{\mathcal{H}}^+(x)$ is connected, $\alpha$ must be surjective. Equation 3 shows that $\alpha^{-1}(y_1)$ is discrete for all $y_1 \in N_{\mathcal{H}}^+(y)$, and with Equation 4, it is elementary to show that $\alpha$ is actually a covering map. We now claim that $N_{\mathcal{H}}^+(y)$ is simply connected. To see this, note that for appropriate $n$, $b_0^n$ maps the ball of any radius in $N_{\mathcal{H}}^+(y)$ into a ball of arbitrarily small radius in $N_{\mathcal{H}}^+(b_0^n(y))$. It follows that $N_{\mathcal{H}}^+(y)$ is a monotone union of open cells. That $N_{\mathcal{H}}^+(y)$ is simply connected now follows from [2].

Since $N_{\mathcal{H}}^+(y)$ is simply connected, it follows that $\alpha$ is a homeomorphism, and therefore invertible. Equations 3 and 4 together now yield

(5) $\qquad d(y_1, y_2) < \zeta$ implies $d(y_1, y_2) \geq (1-\delta) d(\alpha^{-1} y_1, \alpha^{-1} y_2)$.

Using the triangle inequality, we can obtain Equation 5 for all $y_1, y_2 \in N_{\mathcal{H}}^+(y)$; i.e., $\alpha^{-1}$ is a Lipschitz map with Lipschitz constant $\leq \frac{1}{1-\delta}$. As $\delta > 0$ can be chosen arbitrarily small, we conclude $\alpha$ is an isometry. $\square$

*Proof of Theorem* 2.8. For almost every $x \in M$, there exists $a \in \mathcal{H}$ such that the $a$ orbit of $x$ is dense. We will first prove the result for such an $x$, and then complete the proof for arbitrary points in $M$. Thus, we assume that $a \in \mathcal{H}$ and $x \in M$ are chosen so that the $a$ orbit of $x$ is dense. Then for every $y \in N_{\mathcal{H}}^+(x)$, there exists some sequence $\{n_k\}$ such that $a^{n_k}(x) \to y$. By Proposition 2.9, there exists some isometry $\alpha$ of $N_{\mathcal{H}}^+(x)$ such that $\alpha(x) = y$, i.e., $\mathrm{Isom}(N_{\mathcal{H}}^+(x))$ is transitive on $N_{\mathcal{H}}^+(x)$. Since $N_{\mathcal{H}}^+(x)$ is finite dimensional, locally compact, connected and locally connected, this transitive group of isometries is a Lie group [21]. Hence, there exists a $C^\infty$ differentiable structure on $N_{\mathcal{H}}^+(x)$ as a homogeneous space. By [3], these isometries are actually $C^1$ with respect to the original differentiable structure, and by [21, §5.1], it follows that the $C^\infty$ differentiable structure on $N_{\mathcal{H}}^+(x)$ as a homogeneous space is $C^1$ equivalent to the original differentiable structure. Let $g(y, v; t)$ denote the geodesic (with respect to the $C^\infty$ differentiable structure on $N_{\mathcal{H}}^+(x)$ as a homogeneous space) through $y$ with initial velocity $v$. Since the homogeneous metric on each leaf is the restriction of a Hölder metric on all of $M$, it follows that $g$ varies continuously in $y$. Consequently, there exists some $\iota > 0$ such that if $\langle v, v \rangle < \iota$, then for all $y \in M$, $g(y, v; t)$ is defined and is the unique length-minimizing geodesic for all $|t| < 1$.

Let $\mathcal{L}_a(x, y)$ be the set of isometries from $N_{\mathcal{H}}^+(x)$ to $N_{\mathcal{H}}^+(y)$ which can be written as a limit of $\{a^{n_k}\}$ for some sequence $\{n_k\}$. Then $\mathcal{L}_a(x, x)$ is transitive on $N_{\mathcal{H}}^+(x)$. We wish to show first that $\mathcal{L}_a(x,x) \subset \hat{\mathcal{I}}(N_{\mathcal{H}}^+(x))$. To do this, suppose that $\alpha \in \mathcal{L}_a(x,x)$ and $\alpha = \lim_{k\to\infty} a^{n_k}$. Let $x, y \in N_{\mathcal{H}}^+(x)$ such that $d(x,y) < \iota$. Then there exists a unique $v \in \mathcal{N}_{\mathcal{H}}^+(x)$ such that $g(x,v;t)$ is the unique length-minimizing curve from $x$ to $y$. Let $v_k = da^{n_k}(v)$, so that $a^{n_k}(g(x,v;t)) = g(a^{n_k}x, v_k, t)$. Since $\alpha(g(x,v;t))$ is a length-minimizing curve from $\alpha(x)$ to $\alpha(y)$, it follows that there exists $w \in \mathcal{N}_{\mathcal{H}}^+(\alpha(x))$ such that $\alpha(g(x,v;t)) = g(\alpha(x), w; t)$. Since $g$ varies continuously in $x$, we must have $\lim_{k\to\infty} v_k = w$. Since the derivative of an isometry is determined by how geodesics get mapped, we conclude that $d\alpha = \lim_{k\to\infty} da^{n_k}$. Since $a$ preserves $T\mathcal{F}_i$ for every $i$, and $T\mathcal{F}_i$ varies continuously, it follows that $\alpha$ does as well. In other words, $\mathcal{L}_a(x,x) \subset \hat{\mathcal{I}}(N_{\mathcal{H}}^+(x))$. Note that $\hat{\mathcal{I}}(N_{\mathcal{H}}^+(x))$ is a closed subgroup of $\text{Isom}(N_{\mathcal{H}}^+(x))$ since the elements of $\text{Isom}(N_{\mathcal{H}}^+(x))$ are $C^1$. Hence $\hat{\mathcal{I}}(N_{\mathcal{H}}^+(x))$ is a Lie group, and acts transitively on $N_{\mathcal{H}}^+(x)$. Since $N_{\mathcal{H}}^+(x)$ is connected, the connected component $\mathcal{I}(N_{\mathcal{H}}^+(x))$ also acts transitively on $N_{\mathcal{H}}^+(x)$. Indeed, the orbits of $\mathcal{I}(N_{\mathcal{H}}^+(x))$ are open in $N_{\mathcal{H}}^+(x)$, hence closed and by connectedness equal $N_{\mathcal{H}}^+(x)$.

Suppose that $\phi \in \mathcal{I}(N_{\mathcal{H}}^+(x))$ fixes $x$. Since $\phi$ preserves each $T\mathcal{F}_i$, all of which are one dimensional, and $\mathcal{I}(N_{\mathcal{H}}^+(x))$ is connected, $d\phi$ must be the identity. Hence, $\phi$ is the identity and $\mathcal{I}(N_{\mathcal{H}}^+(x))$ acts without isotropy on $N_{\mathcal{H}}^+(x)$.

To complete the proof consider an arbitrary $z \in M$. As in the proof of Proposition 2.9, the density of the $a$ orbit through $x$ implies there exists an isometry of $N_{\mathcal{H}}^+(x)$ with $N_{\mathcal{H}}^+(z)$; i.e. $\mathcal{L}(x,z)$ is nonempty. Using the argument above, we may conclude that for any $\theta \in \mathcal{L}(x,z)$ and $\phi \in \mathcal{L}(x,x)$, $\theta\phi\theta^{-1} \in \mathcal{I}(N_{\mathcal{H}}^+(z))$. The result now follows. □

*Remark* 2.13. Later we shall need to make use of the fact that this argument applies when $M$ is the fiber of some bundle $X \to B$. More specifically, assume that $A$ acts via bundle automorphisms on the bundle $X \to B$ with compact fibers $M$ satisfying Assumptions (A2) and (A3). If, in addition, the action is trellised with respect to $A$ in the direction of the fibers and there exists the appropriate $A$ equivariant Hölder Riemannian metric on each fiber, then our argument still holds. The key ingredients are density of $a$ orbits in $X$ and compactness of the fiber $M$. Compactness of the base space is irrelevant.

2.3. *Unions of stable and unstable foliations.* Having identified the structure of certain types of stable submanifolds, we will now attempt to do the same for a class of more general sets. To begin, we present some necessary technical facts.

LEMMA 2.14. *Let $c \in A$ be regular, let $d_c^-$ be the leafwise distance for the $W_c^-$ foliation, and define $W_{c,\varepsilon}^-(x) = \{y \in W_c^-(x) | d_c^-(x,y) < \varepsilon\}$.*

1. *There exists $\varepsilon > 0$ such that for every $x, y \in M$, the intersection of $W^{0+}_{c,\varepsilon}(y) \cap W^-_{c,\varepsilon}(x)$ consists of at most a single point, called $[x, y]_c$.*

2. *There exists $\delta > 0$ such that if $d(x, y) < \delta$, then $W^{0+}_{c,\varepsilon}(y) \cap W^-_{c,\varepsilon}(x) = \{[x, y]_c\}$.*

3. *There exists $\delta_c > 0$ such that if $d(x, y) < \delta_c$, then $d(x, [x, y]_c), d(y, [x, y]_c) < \frac{\varepsilon}{m(c)}$, where $m(c)$ is the maximal expansion of $dc$ on $TM$ with respect to the Hölder metric.*

*Proof.* See Proposition 6.4.13 in [14] and Theorem 6.1.9 in [9]. □

LEMMA 2.15. *If $c_1, c_2 \in A$ are two regular elements, and $a$ is a positive linear combination of $c_1$ and $c_2$, then there exists $\varepsilon > 0$ such that for every $z \in W^-_{c_1,\varepsilon}(x) \cap W^-_{c_2,\varepsilon}(x)$, there exist constants $C > 0$ and $\lambda > 0$ so that $d(a^n z, a^n x) < C e^{-n\lambda} d(z, x)$.*

*Proof.* Without loss of generality we may assume $a = c_1 + c_2$ by passing to suitable roots of $c_1$ and $c_2$ in $A$. Since the $c_i$ are regular, there exist constants $\varepsilon > 0$, $C_i > 0$ and $\lambda_i > 0$ such that for every $x \in M$ and for every $y \in W^-_{c_i,\varepsilon}(x)$, $d(c_i^n y, c_i^n x) < C_i e^{-n\lambda_i} d(y, x)$. Since $c_1$ and $c_2$ commute, $c_1 W^-_{c_2}(x) = W^-_{c_2}(c_1^n x)$. Hence, letting $C = \max(C_1, C_2)$ and $\lambda = \min(\lambda_1, \lambda_2)$ the lemma follows. □

Given an $a \in \mathcal{H}$, there exist normally hyperbolic $a_1, a_2$ nearby such that $E^-_a(x) = E^-_{a_1}(x) \cap E^-_{a_2}(x)$ and $E^+_a(x) = E^+_{a_1}(x) \cap E^+_{a_2}(x)$, and $\mathcal{N}^{\pm}_{\mathcal{H}}(x)$ is the neutral direction (cf. proof of Lemma 2.7); i.e., $\mathcal{N}^-_{\mathcal{H}}(x) = W^-_{a_1}(x) \cap W^+_{a_2}(x) \cap W^-_{b_0}(x) = W^-_{a_1}(x) \cap W^+_{a_2}(x)$. It is also clear that $E^-_{a_i}(x)$ is tangent to the foliation $W^-_{a_i}(x)$ and similarly for $E^+_{a_i}(x)$. As a result, $E^-_a(x)$ is tangent to the foliation $W^-_{a_1}(x) \cap W^-_{a_2}(x)$, and similarly for the unstable directions. In particular, $E^+_a(x)$ is integrable with integral foliation $W^+_a(x)$ defined to be $W^+_{a_1}(x) \cap W^+_{a_2}(x)$. Note that by Lemma 2.15, $W^+_a(x)$ is indeed contracted by $a$, because both $a_1$ and $a_2$ contract $W^+_a(x)$ and $a$ can be written as a positive linear combination of $a_1$ and $a_2$.

LEMMA 2.16. *Let $c \in A$ be regular, and choose $\delta_{b_0}$ as in Lemma 2.14.3. Suppose $y \in W^-_{b_0}(x)$ with $d(x, y) < \delta_{b_0}$. Then $[x, y]_c \in W^-_{b_0}(x)$.*

*Proof.* First we claim that $b_0[x, y]_c = [b_0 x, b_0 y]_c$ if $d(x, y) < \delta_{b_0}$. Clearly $b_0[x, y]_c \in W^-_c(x) \cap W^{0+}_c(y)$, and by choice of $\delta_{b_0}$, we have $d(b_0 x, [x, y]_c), d(b_0 y, [x, y]_c) < \frac{\varepsilon}{m(b_0)} m(b_0) = \varepsilon$. So, by uniqueness, our claim follows.

Note that $y \in W^-_{b_0}(x)$ implies $d^-_c(b_0 x, b_0 y) < d^-_c(x, y)$. Since locally we can bound $d(b_0 x, b_0 y)$ in terms of $d^-_c(b_0 x, b_0 y)$, we have that $d(b_0 x, b_0 y) < d(x, y)$.

This allows us to repeat the argument, yielding $b_0^n[x,y]_c = [b_0^n x, b_0^n y]_c$. As $b_0^n x$ and $b_0^n y$ approach each other, it follows that $[b_0^n x, b_0^n y]_c$ approaches $b_0^n x$; i.e., $[x,y]_c \in W_{b_0}^-(x)$. □

LEMMA 2.17. *Suppose $a \in \mathcal{H}$. There exists $\Delta > 0$ such that if $d(a^n x, a^n y) < \Delta$ for all $n \in \mathbb{Z}$, then $[x,y]_{b_0} \in N_{\mathcal{H}}^-(x)$ and $y \in AN_{\mathcal{H}}^+([x,y]_{b_0})$.*

*Proof.* Recall that $N_{\mathcal{H}}^-(x) = W_{a_1}^-(x) \cap W_{a_2}^+(x) \cap W_{b_0}^-(x)$. Without loss of generality we may assume that $a_1$ and $a_2$ are chosen so that $a$ is a positive linear combination of $a_2$ and $b_0$, but not a positive linear combination of $a_1$ and $b_0$. Let $m = \max(m(b_0), m(a_1), m(a_2))$. Choose $\varepsilon$ so that $\varepsilon/m < \min(\delta_{b_0}, \delta_{a_1}, \delta_{a_2})$. For $w$ sufficiently close to $z$, define

$$\begin{aligned}\mathcal{S}(w) &= \max\{d(w, [w,z]_{b_0}), d(w, [w,z]_{a_1}), d(w, [w,z]_{a_2}), \\ &\quad d(z, [w,z]_{b_0}), d(z, [w,z]_{a_1}), d(z, [w,z]_{a_2})\}.\end{aligned}$$

By Lemma 2.14, we can pick $\Delta_1 < \varepsilon/m$ so that if $d(w,z) < \Delta_1$, then $\mathcal{S}(w) < \varepsilon/m$. Next, pick $\Delta < \varepsilon/m$ so that if $d(w,z) < \Delta$, then $\mathcal{S}(w) < \Delta_1$.

Let $z = [x,y]_{b_0}$. First, we claim that $z \in W_{a_1}^-(x)$. In particular, we will show that if $z_1 = [x,z]_{a_1}$, then $z = z_1$. By Lemma 2.16, $z_1 \in W_{b_0}^-(x) = W_{b_0}^-(z)$. However, by definition, $z_1 \in W_{a_1}^{0+}(z)$. Since $b_0$ is regular, we may assume that $z_1 \in W_{a_1}^+(z)$. Consequently, by Lemma 2.15, $z_1 \in W_a^+(z)$. Now, since $\{d(a^n x, a^n y)\}$ is bounded by $\Delta$, $\{d(a^n x, a^n z)\}$ must be bounded by $\Delta_1 < \frac{\varepsilon}{m}$, and we may conclude that $[a^n x, a^n z]_{a_1}$ is defined for all $n$, and that $\{d([a^n x, a^n z]_{a_1}, a^n z)\}$ is bounded by $\frac{\varepsilon}{m}$. By the uniqueness of canonical coordinates, we must have $[a^n x, a^n z]_{a_1} = a^n z_1$. Thus $\{d(a^n z_1, a^n z)\}$ is bounded by $\frac{\varepsilon}{m}$ for all $n \in \mathbb{Z}$.

As $z_1 \in W_a^+(z)$, the leafwise distance between $z_1$ and $z$ grows under positive iterates of $a$, unless $z = z_1$. Since the foliations are Hölder, there exists a neighborhood of $x$ in the leaf in which the ambient distance is a continuous function of the leafwise distance. By shrinking $\varepsilon$ (and hence $\Delta$) if necessary, we may assume that the $\varepsilon/m$-ball about $x$ lies in this neighborhood. Thus, if $z_1 \in W_a^+(z)$ then there exists some $n$ so that $d(a^n z, a^n z_1) > \varepsilon/m$ unless $z = z_1$. We are forced to conclude that $z = z_1$.

Now let $z_2 = [x,z]_{a_2}$. By definition, $z_2 \in W_{a_2}^-(x)$, and since $z \in W_{b_0}^-(x)$, Lemma 2.16 implies $z_2 \in W_{b_0}^-(x)$. Thus $z_2 \in W_a^-(x)$. As above, we see that $a^n z_2 = [a^n x, a^n z]_{a_2}$, and that $\{d(a^n z_2, a^n x)\}$ is bounded by $\frac{\varepsilon}{m}$ for all $n \in \mathbb{Z}$. However, since $z_2 \in W_a^-(x)$, the distance between $x$ and $z_2$ will be expanded under negative iterates of $a$, unless $x = z_2$. Again, this contradiction forces $x = z_2$. We conclude $z \in W_{a_1}^-(x) \cap W_{a_2}^+(x) \cap W_{b_0}^-(x) = N_{\mathcal{H}}^-(x)$.

The second claim follows in a similar manner. □

Define $P_\mathcal{H}$ to be the set of paths in $M$ which are piecewise tangent to $N_\mathcal{H}^+$ or $A \cdot N_\mathcal{H}^-$. More explicitly, if $\gamma : [0,1] \to M$ is a path, then $\gamma \in P_\mathcal{H}$ if there exists a sequence $0 = t_0 < t_1 < \cdots < t_k = 1$ such that $\gamma|_{[t_i, t_{i+1}]}$ is $C^1$, $\gamma|_{[t_i, t_{i+1}]} \subset N_\mathcal{H}^+(t_i)$ or $\gamma|_{[t_i, t_{i+1}]} \subset A \cdot N_\mathcal{H}^-(t_i)$. Next, we define

$$N_\mathcal{H}(x) = \{y \in M | \text{ there exists } \gamma \in P_\mathcal{H} \text{ with } \gamma(0) = x \text{ and } \gamma(1) = y\}.$$

We intend to demonstrate the existence of a smooth differentiable structure on $N_\mathcal{H}(x)$. The first step in this direction is to define a distance function on $N_\mathcal{H}(x)$. For $x_1, x_2 \in N_\mathcal{H}(x)$ define $\hat{d}(x_1, x_2)$ to be the infimum of the lengths measured via the Hölder metric of any path $\gamma \in P_\mathcal{H}$ with $\gamma(0) = x_1$ and $\gamma(1) = x_2$. Note that $\hat{d}(x_1, x_2) \geq d(x_1, x_2)$. It is then a simple exercise to prove:

LEMMA 2.18. $\hat{d}$ *defines a metric on* $N_\mathcal{H}(x)$.

PROPOSITION 2.19. *The topology generated by $\hat{d}$ makes $N_\mathcal{H}(x)$ a connected, locally connected, locally compact, finite dimensional, locally simply connected topological space.*

*Proof.* Note that any $\gamma \in P_\mathcal{H}$ is continuous under $\hat{d}$. Connectedness and local connectedness thus follow.

To see that the topology is locally compact, consider a sequence $\{x_n\} \subset N_\mathcal{H}(x)$ such that $\hat{d}(x_n, x) \leq \zeta$ for some sufficiently small $\zeta$. Since $\hat{d} \geq d$, we may assume without loss of generality that $x_n$ converges to some $y$ with respect to $d$. To establish local compactness, we first show that $y \in N_\mathcal{H}(x)$ and second show that $\hat{d}(x_n, y) \to 0$ as $n \to \infty$. Then $\hat{d}(x, y) \leq \zeta$, and the $\zeta$-ball about $x$ is compact.

If $\zeta$ is chosen sufficiently small, then Lemma 2.17 implies that $[x, x_n]_{b_0} \in N_\mathcal{H}^+(x) \cap A \cdot N_\mathcal{H}^-(x_n)$. Further, we have that $d([x, x_n]_{b_0}, [x, y]_{b_0}) \to 0$ as $n \to \infty$. Since the closed ball about $x$ in $N_\mathcal{H}^+(x)$ is closed in the Hölder metric, we conclude $[x, y]_{b_0} \in N_\mathcal{H}^+(x)$. Also, since locally the leaves of $A \cdot N_\mathcal{H}^-$ are continuous in a sufficiently small neighborhood, the fact that $[x, x_n]_{b_0} \in A \cdot N_\mathcal{H}^-(x_n)$ implies that $[x, y]_{b_0} \in A \cdot N_\mathcal{H}^-(y)$. Hence $[x, y]_{b_0} \in N_\mathcal{H}^+(x) \cap A \cdot N_\mathcal{H}^-(y)$ from which we conclude $y \in N_\mathcal{H}(x)$.

Let $\delta_{b_0}$ be as in Lemma 2.14.3. Since $\{x_n\} \subset N_\mathcal{H}(x)$, Lemma 2.17 implies there exists some $\xi$ such that if $d(x_k, x) < \xi$, then $d(a^n x, a^n x_k) < \delta_{b_0}/2$ for all $n \in \mathbb{Z}$. Similarly, since we may pick $\zeta < \xi/2$, $d(a^n x, a^n y) < \delta_{b_0}/2$ for all $n \in \mathbb{Z}$. It follows that for all $x_k$ sufficiently close to $y$, $d(a^n x_k, a^n y) < \delta_{b_0}$ for all $n \in \mathbb{Z}$. As a consequence of Lemma 2.17, $[x_k, y]_{b_0} \in N_\mathcal{H}^+(x_k) \cap A \cdot N_\mathcal{H}^-(y)$ for $x_k$ sufficiently close to $y$. Since the leaves of the $N_\mathcal{H}^+$ foliation are uniformly Hölder continuous, we have that $\hat{d}(x_k, [x_k, y]_{b_0}) \to 0$ as $n \to \infty$. Similarly, using the continuity of $[\cdot, \cdot]_{b_0}$, we have $\hat{d}([x_k, y]_{b_0}, y) \to 0$ as $n \to \infty$. Thus $\hat{d}(x_k, y) \to 0$ as $n \to \infty$, and local compactness of the topology results.

To see that the topology is finite dimensional, consider a closed ball about $x$ in $N_{\mathcal{H}}(x)$. Because $\hat{d} \geq d$, the canonical embedding of such a ball into $M$ is continuous. Thus, the image of a closed ball of sufficiently small radius about $x$ is homeomorphic onto its image. Since any compact subset of Euclidean space is finite dimensional, it follows that $N_{\mathcal{H}}(x)$ is finite dimensional.

To see that the topology is locally simply connected, note that Lemma 2.17 implies that $N_{\mathcal{H}}(x)$ has a local product structure with respect to $N_{\mathcal{H}}^+(x)$ and $AN_{\mathcal{H}}^-(x)$. Hence, it must be locally simply connected. □

PROPOSITION 2.20. *With respect to the metric $\hat{d}$, there exists a transitive set of local isometries on $N_{\mathcal{H}}(x)$ for almost every $x$. Hence, $N_{\mathcal{H}}(x)$ is a locally homogeneous space for every $x$.*

*Proof.* The proof is essentially the same as for Proposition 2.9 with minor modifications. We begin by considering an $x \in M$ which lies in the conull set of points in $M$ such that there exists $a \in \mathcal{H}$ so that the $a$-orbit of $x$ is dense. First, note that if $\gamma \in P_{\mathcal{H}}$ then $a(\gamma) \in P_{\mathcal{H}}$, and therefore $\alpha(\gamma) \in P_{\mathcal{H}}$. The ergodicity of $a \in A$ implies that for every $y \in N_{\mathcal{H}}(x)$, there exists an $\alpha$ such that $\alpha(x) = y$. Second, since $N_{\mathcal{H}}(x)$ need not be simply connected, the proof of Proposition 2.9 implies only that $\alpha$ is a covering map and not necessarily a homeomorphism. However, the rest of the proof shows that $\alpha$ is a local isometry. By lifting $\alpha$ to $\tilde{N}_{\mathcal{H}}(x)$, the universal cover of $N_{\mathcal{H}}(x)$, we obtain a map $\tilde{\alpha}$. The arguments of Proposition 2.9 now imply that $\tilde{\alpha}$ is an isometry with respect to the distance function lifted from $N_{\mathcal{H}}(x)$. Since the deck transformations are also isometries, we conclude that $\tilde{N}_{\mathcal{H}}(x)$ admits a transitive group of isometries.

That we can write $\tilde{N}_{\mathcal{H}}(x) \cong L/L_x$ follows from the work of Montgomery and Zippin [21]. Consequently, $N_{\mathcal{H}}(x)$ is a locally homogeneous space.

For an arbitrary point $y \in M$, we find a sequence $\{n_k\}$ so that for $x$ as above, we have $a^{n_k}x \to y$. Using the same limiting argument from the proof of Proposition 2.9, we obtain a local isometry $\beta : N_{\mathcal{H}}(x) \to N_{\mathcal{H}}(y)$ which lifts to an isometry $\tilde{\beta} : \tilde{N}_{\mathcal{H}}(x) \to \tilde{N}_{\mathcal{H}}(y)$. Hence, we can write $\tilde{N}_{\mathcal{H}}(y)$ as a homogeneous space, and $N_{\mathcal{H}}(y)$ as a locally homogeneous space. □

For all $y \in N_{\mathcal{H}}(x)$, we have $N_{\mathcal{H}}^+(y) \subset N_{\mathcal{H}}(x)$. The following proposition shows that this inclusion is well behaved.

PROPOSITION 2.21. *For every $y \in N_{\mathcal{H}}(x)$, $N_{\mathcal{H}}^+(y)$ is a $C^\infty$ submanifold of $N_{\mathcal{H}}(x)$ (with respect to its differentiable structure as a homogeneous space) which gives rise to an L-invariant foliation.*

*Proof.* Let $i_y : N_{\mathcal{H}}^+(y) \hookrightarrow N_{\mathcal{H}}(x)$ be the canonical inclusion. Recall that the leaves of $N_{\mathcal{H}}^+$ admit distance functions defined by restriction of the ambient Hölder metric, and $N_{\mathcal{H}}(x)$ has the distance function $\hat{d}$ defined above. With respect to these distances, $i_y$ is distance nonincreasing, thus Lipschitz, and therefore almost everywhere differentiable. Using transitivity of $\mathcal{I}(N_{\mathcal{H}}^+(y))$ on $N_{\mathcal{H}}^+(y)$, for every $z \in N_{\mathcal{H}}^+(y)$ there exists a distribution tangent to $N_{\mathcal{H}}^+$. Since this distribution is invariant under the action of $L$, it must be $C^\infty$. The proposition now follows. □

*Remark* 2.22. We interpret Propositions 2.20 and 2.21 as providing a smooth differentiable structure on $N_{\mathcal{H}}(x)$ with respect to which the Lyapunov spaces lying in $N_{\mathcal{H}}(x)$ vary smoothly with respect to each other. Of course, it remains to be seen whether this new differentiable structure has any relation to the original differentiable structure on $M$.

2.4. *Smooth conjugacy.* The next goal in our exposition is to show that the differentiable structure on $N_{\mathcal{H}}(x)$ as a locally homogeneous space is $C^\infty$ equivalent to the original differentiable structure on $M$. To do this we will adapt an argument by Katok and Lewis [16]. The first step in this argument is to use the nonstationary Sternberg linearization to construct a smooth conjugacy between the one dimensional foliations $\mathcal{F}_i$.

Let $\mathcal{L} = M \times \mathbb{R}$ be the trivial line bundle over $M$. Fix some smooth Riemannian metric on $M$, and define $\mathcal{E}_i : \mathcal{L} \to M$ so that $\mathcal{E}_i(x,t)$ is the point on $\mathcal{F}_i(x)$ which is $t$ units from $x$ measured with respect to this smooth Riemannian metric. We may have to pass to a finite cover of $M$ to ensure that there exists an orientation for $\mathcal{F}_i$, for all $i$. We note that $\mathcal{E}_i$ varies smoothly with $t$, and by Assumption (A4), varies Hölder with $x$.

For any $a \in A$ define
$$\hat{a} : \mathcal{L} \to \mathcal{L} \text{ so that } \hat{a}(x,t) = (ax, e^{\chi_i(a)}t),$$
and also define
$$\tilde{a} : \mathcal{L} \to \mathcal{L} \text{ so that } \tilde{a}(x,t) = (ax, \mathcal{E}_i^{-1}(ax, a\mathcal{E}_i(x,t))),$$
where $\mathcal{E}_i^{-1}$ is the leafwise inverse. Next define $H : \mathcal{L} \to \mathcal{L}$ so that $H(x,t) = (x, d_H(x, \mathcal{E}_i(x,t)))$ where $d_H$ is the distance measured with respect to the Hölder metric.

By the nonstationary Sternberg linearization [16], there exists a unique reparametrization $G : \mathcal{L} \to \mathcal{L}, (x,t) \mapsto (x, G_x(t))$ such that

1. each $G_x$ is a $C^\infty$ diffeomorphism of $\mathbb{R}$,

2. $G_x(0) = 0, G'_x(0) = 1$ for every $x \in M$,

3. $x \mapsto G_x$ is a continuous map $M \to C^\infty(\mathbb{R})$, and

4. $G\tilde{a}(x,t) = \hat{a}G(x,t)$ for every $x \in M, t \in \mathbb{R}$.

As in [16, Lemma 4.8], we remark that since $A$ is abelian, the uniqueness property of $G$ ensures that $G$ simultaneously linearizes the transformations on $\mathcal{L}$ for every $a \in A$. We claim that $H$ satisfies the same equivariance property as $G$, i.e., $H\tilde{a} = \hat{a}H$.

$$\begin{aligned} H\tilde{a}(x,t) &= H(ax, \mathcal{E}_i^{-1}(ax, a\mathcal{E}_i(x,t))) \\ &= (ax, d_H(ax, \mathcal{E}_i(ax, \mathcal{E}_i^{-1}(ax, a\mathcal{E}_i(x,t))))) \\ &= (ax, d_H(ax, a\mathcal{E}_i(x,t))) = (ax, e^{\chi_i(a)}d_H(x, \mathcal{E}_i(x,t))) \\ &= \hat{a}(x, d_H(x, \mathcal{E}_i(x,t))) = \hat{a}H(x,t). \end{aligned}$$

Let $K = GH^{-1}$. Then we have

$$K\hat{a} = GH^{-1}\hat{a} = G\tilde{a}H^{-1} = \hat{a}GH^{-1} = \hat{a}K.$$

Hence, if $K(x,t) = (x, K_x(t))$, then we have

$$e^{\chi_i(a)}K_x(t) = K_{ax}(e^{\chi_i(a)}t).$$

However, there exists some $j$ such that $T\mathcal{F}_i \subset E_j$. So, if we choose $a \in \ker(\chi_j)$, then $K_x(t) = K_{ax}(t)$. Note that the map $\kappa : M \to C^{\text{Hölder}}(\mathbb{R}), \kappa(x) = K_x$, is a continuous map into a countably separated space. By Assumption (A3), the action of the one-parameter subgroup of $A$ through $a$ is ergodic. Hence $\kappa$ must be constant; i.e., for every $x \in M$, $K_x(t) = K(t)$ for some Hölder map $K : \mathbb{R} \to \mathbb{R}$. It follows that $K(e^{\chi_j(a)}t) = e^{\chi_j(a)}K(t)$ for all $t \in \mathbb{R}$ and any $a \in A$. If $a$ is regular, then it is easy to see that $K$ must be $C^\infty$ for all $t$ except perhaps 0. However, smoothness at 0 follows exactly as in [15]. Since $H_x^{-1} = G_x^{-1}K$ determines the inclusion of $i_\mathcal{H}(x)$ along $\mathcal{F}_i$, we have just proved:

PROPOSITION 2.23. *Suppose $T\mathcal{F}_i \subset E_j, \mathcal{H} \subset \ker(\chi_j)$, and $i_\mathcal{H}(x): N_\mathcal{H}(x) \to M$ is the inclusion map. Then*

1. $i_\mathcal{H}(x)$ *is $C^\infty$ along $\mathcal{F}_i(y)$ for every $y \in N_\mathcal{H}(x)$, and*

2. *the map $y \mapsto i_\mathcal{H}(x)|_{\mathcal{F}_i(y)}$ is uniformly continuous in the $C^\infty$ topology on $C^\infty(\mathbb{R}, M)$ where $y$ ranges over $N_\mathcal{H}(x)$.*

By a similar argument as for $K$ we see that $K^{-1}$ is smooth, which implies that $K$ is a diffeomorphism. Hence we get the following:

COROLLARY 2.24. *Let $\mathcal{H} \subset A$ be a linear subspace. Then $i_\mathcal{H}(x) : N_\mathcal{H}(x) \to M$ is a $C^\infty$ immersion for every $x \in M$.*

The proof of this corollary is similar to an argument due to Katok and Lewis [16], which requires the following result of Journé [13].

LEMMA 2.25. *Let $M$ be a $C^\infty$ manifold and $\mathcal{L}$ and $\mathcal{L}'$ two Hölder foliations, transverse, and with uniformly $C^\infty$ leaves. If a function $f$ is uniformly $C^\infty$ along the leaves of the two foliations, then it is $C^\infty$ on $M$.*

*Proof of Corollary* 2.24. Using Remark 2.6, there exists a nested sequence of Hölder foliations $\mathcal{L}_1 \subset \mathcal{L}_2 \subset \cdots \subset N_\mathcal{H}^+(x)$ with uniformly $C^\infty$ leaves such that $\bigoplus_{j=1}^i E_{b_j}^+ \cap \mathcal{N}_\mathcal{H}^+$ is the distribution tangent to $\mathcal{L}_i$. Using the arguments following Lemma 4.11 in [16], we can apply Lemma 2.25 inductively to see that $i_\mathcal{H}(x)$ is smooth along each leaf of $\mathcal{L}_i$ for every $i$, and therefore is $C^\infty$ along $N_\mathcal{H}^+(y)$ for every $y \in N_\mathcal{H}(x)$. An analogous argument allows us to conclude that $i_\mathcal{H}(x)$ is $C^\infty$ along $N_\mathcal{H}^-(y)$ for every $y \in N_\mathcal{H}(x)$. Next, since $i_\mathcal{H}(x)$ must be $C^\infty$ along the $A$ orbit, Lemma 2.25 implies that $i_\mathcal{H}(x)$ is smooth along the weak unstable leaves in $N_\mathcal{H}(x)$. A final application of Lemma 2.25 to the weak unstable and stable foliations shows that $i_\mathcal{H}(x)$ is $C^\infty$ on all of $N_\mathcal{H}(x)$. □

*Proof of Theorem* 2.4. The proof is essentially the same as for Corollary 2.24. Let $\chi$ and $\lambda$ be two distinct weights for $A$. Since $k \geq 3$, there exists a nontrivial linear subspace $\mathcal{H} \subset A$ such that $\mathcal{H} \subset \ker(\chi) \cap \ker(\lambda)$. Thus, the distributions $E_\chi$ and $E_\lambda$ both lie in the distributions tangent to $N_\mathcal{H}$. Corollary 2.24 implies that $E_\chi$ varies smoothly with respect to the foliation corresponding to $E_\lambda$ for any other weight $\lambda$. Using Lemma 2.5 and applying Lemma 2.25 inductively, it follows that the distribution $E_\chi$ varies smoothly along the stable foliation for $b_0$. Again, a similar argument can be used to show that $E_\chi$ varies smoothly along the weak unstable foliation for $b_0$, and, with yet another application of Lemma 2.25, we have that $E_\chi$ varies smoothly along $M$.

To see that the Hölder metric is actually $C^\infty$, choose a section of the frame bundle which is orthonormal with respect to this metric at every point. By Proposition 2.23.1 and the fact that the $E_\chi$ now all vary smoothly, it follows that this section, and hence the metric, must be smooth. □

## 3. Actions of semisimple groups and their lattices

We now turn our attention to the primary subjects of this paper, Anosov actions of semisimple groups and their lattices on closed manifolds. By considering actions of abelian subgroups we can apply results from Section 2 to prove Theorems 1.2 and 1.5.

3.1. *Trellised actions of semisimple Lie groups.* Let $G$ be a connected semisimple Lie group without compact factors with real rank at least three, and let $A \subset G$ be a maximal $\mathbb{R}$-split Cartan subgroup. The following two results, taken from [7], allow us to apply our results from Section 2.

THEOREM 3.1. *Suppose $G$ is a connected semisimple Lie group of higher rank without compact factors such that each simple factor of $G$ has $\mathbb{R}$-rank at least $2$. Suppose that $G$ acts on a closed manifold $M$ such that the $G$ action is Anosov and volume-preserving. Let $H$ be the Hölder algebraic hull of the $G$ action on $P \to M$, the $G$-invariant reduction of the derivative action on the full frame bundle over $M$. Then, by possibly passing to a finite cover of $G$, there exist*

1. *a normal subgroup $K \subset H$,*

2. *a Hölder section $s : M \to P/K$, and*

3. *a homomorphism $\pi : G \to H/K$,*

*such that $s(gm) = g.s(m).\pi(g)^{-1}$ for every $g \in G$ and every $m \in M$.*

*Moreover, if the irreducible subrepresentations of $\pi$ are multiplicity free, then $K \subset H$ is a compact normal subgroup.*

*Remark* 3.2. Note that the Hölder algebraic hull $H$ may be viewed as a subgroup of $\mathrm{SL}(n, \mathbb{R})$. By [7, Th. 3.1], $H$ is reductive with compact center. Using this, and the fact that $K$ is normal, we can produce a homomorphism (which by abuse of notation we also call) $\pi : G \to H \subset \mathrm{SL}(n, \mathbb{R})$ which modulo $K$ factors to $\pi$. This is the representation to which we refer in the above theorem. See [7, Rem. 3.5].

*Definition* 3.3. We call a volume-preserving Anosov action *multiplicity free* if the irreducible subrepresentations of $\pi$ are multiplicity free.

COROLLARY 3.4. *Let $G$, $P$, $M$, and $H$ be as in Theorem* 3.1. *Assume that the action is multiplicity free so that $K$ is compact. Let $A$ be a maximal $\mathbb{R}$-split Cartan of $G$, with $\{\chi\}$ the set of weights of $\pi$ with respect to $A$. There exist*

1. *a $K$-invariant Hölder Riemannian metric, $\|\cdot\|_K$, on $M$, and*

2. *a $K$-invariant Hölder decomposition $TM = \bigoplus E_\chi$*

*such that for every $v \in E_\chi$ and $a \in A$*

$$\|av\|_K = e^{\chi(\log a)} \|v\|_K.$$

Combining Corollary 3.4 with Theorem 2.4 we obtain:

COROLLARY 3.5. *Suppose $G$ is a connected semisimple Lie group without compact factors such that the $\mathbb{R}$-rank of $G$ is at least $3$. Suppose that $G$ acts on a closed manifold $M$ such that the $G$ action is Anosov and volume-preserving. Let $A$ be a maximal $\mathbb{R}$-split Cartan of $G$, with $\{\chi_i\}$ the set of weights of $\pi$ with* respect to $A$. If the action of $G$ on $M$ is multiplicity free and trellised with respect to $A$, then there exist

1. a $C^\infty$ Riemannian metric, $\|\cdot\|$, on $M$, and

2. a $C^\infty$ decomposition $TM = TA \bigoplus E_i$, where $TA$ is the tangent space to the $A$ orbit,

such that for every $v \in E_i$ and $a \in A$

$$\|av\| = e^{\chi_i(\log a)}\|v\|.$$

Moreover, if the trellis consists of one dimensional strongest stable foliations, i.e., the action is Cartan, then this result holds when the real rank of $G$ is at least two.

*Proof.* The result follows immediately by application of Theorem 2.4; hence we need only check that Assumptions (A0) through (A6) from Section 2 hold. Any Anosov action satisfies (A6). Assumption (A0) holds since Anosov actions are locally faithful by definition. Assumption (A2) holds by hypothesis. Since $A \subset G$ is a noncompact subgroup, Assumption (A3) follows from Moore's Ergodicity Theorem. Finally, Corollary 3.4 ensures that Assumptions (A1), (A4) and (A5) hold.

For the Cartan case, we use the smoothness of the one dimensional strongest stable foliations inside a stable manifold (cf. proof of Theorem 3.12). □

COROLLARY 3.6. *Assume the conditions of Corollary* 3.5. *Then, by possibly passing to a finite cover of* $M$, *there exists a* $C^\infty$ *totally* $\pi$-*simple framing of* $M$.

*More specifically, let* $P \to M$ *be the principal* $H$ *bundle which is the* $G$-*invariant reduction of the derivative action on the full frame bundle over* $M$ *as in Theorem* 3.1. *If* $M' \to M$ *is the finite cover, and if* $P' \to M'$ *is the principal* $H$ *bundle over* $M'$ *lifted from* $P$, *then there exists a* $C^\infty$ *section* $\phi$ *of* $P' \to M'$ *such that* $\phi(gm) = g.\phi(m).\pi(g)^{-1}$ *for all* $g \in G$ *and for all* $m \in M'$.

*Proof.* Let $M' \to M$ be the finite cover such that for each $i$ the foliation $\mathcal{F}_i$ is orientable. For each $i = 1, \ldots, n$, define a vector field $X_i$ on $M'$ such that $X_i(x)$ is an element of $T\mathcal{F}_i(x)$ with unit length measured with respect to the $C^\infty$ Riemannian metric from Corollary 3.5. Smoothness of the foliations of the trellis and the metric assure us that each $X_i$ is indeed $C^\infty$. Hence, we have a $C^\infty$ section $\phi$ of the full frame bundle $P' \to M'$.

Let $p : P' \to P'/K$ be the usual projection, where $K$ is as in Corollary 3.4. From the proofs of Theorem 3.1 and Corollary 3.4 in [7], it follows that $\phi$ projects to a totally $\pi$-simple section $s' : M' \to P'/K$ induced from $s : M \to P/K$, i.e., $p \circ \phi = s'$. It follows that $\phi$ may be constructed from $s'$ along with some smooth cocycle $\kappa : G \times M \to K$. More explicitly, we have $\phi(gx) = dg.\phi(x).\pi(g)^{-1}.\kappa(g,x)$. To show that $\phi$ is totally $\pi$-simple, and thereby complete the proof, we need to demonstrate that $\kappa$ is trivial.

First, note that $\kappa : A \times M \to K$ must be trivial. Without loss of generality, we may assume that $\pi(A) \subset H \subset \mathrm{SL}(n,\mathbb{R})$ is a diagonal subgroup, and, in particular, that each entry along the diagonal of $\pi(a)$ has the form $e^{\chi_l(a)}$ for some weight $\chi_l$. Since $T\mathcal{F}_i \subset E_j$ for some weight $\chi_j$, Assumption (A5) implies that $da(X_i(x)) = e^{\chi_j(a)}.X_i(ax)$. By construction of $\phi$, we have $\phi(ax) = da.\phi(x).\pi(a)^{-1}$, forcing $\kappa(a,x) = 1$ for all $a \in A$ and all $x \in M$.

From the cocycle identity, for any $g \in G$, we have $\kappa(a^{-n}ga^n, x) = \kappa(g, a^n x)$. Suppose that $u \in G$ is unipotent. We can find $a \in A$ so that $a^{-n}ua^n \to 1$ as $n \to \infty$. By continuity of $\kappa$, we have $\kappa(a^{-n}ua^n, x) \to 1$ as $n \to \infty$. Now, if $x \in M$ is a recurrent point, there exists some subsequence $\{n_k\}$ such that $a^{n_k}x \to x$ as $k \to \infty$. Hence, as $k \to \infty$, $\kappa(u, a^{n_k}x) \to \kappa(u,x)$. Since $\kappa(a^{-n}ua^n, x) = \kappa(u,x)$, we conclude $\kappa(u,x) = 1$. But the set of recurrent points in $M$ is dense, so continuity of $\kappa$ forces $\kappa$ to be trivial for any unipotent element. As $G$ is generated by this set, it follows that $\kappa$ is the trivial cocycle. □

*Proof of Theorem* 1.2. The proof is completed by applying arguments from [25, §3] and Feres [4]. As an aid to the reader, we will briefly summarize these arguments here.

From Corollary 3.6, we have a $C^\infty$ framing for $M'$. The first step is to show that this framing, viewed as a collection of nonvanishing $C^\infty$ vector fields, generates a finite dimensional Lie algebra. Let $\{X_i\}$ be the $C^\infty$ vector fields determined by the framing, and also identify $T_xM \cong \mathbb{R}^n$ via this framing in the standard way. Define $f(x) : \mathbb{R}^n \times \mathbb{R}^n \to \mathbb{R}^n$ so that $f(x)((u_1,\ldots,u_n),(v_1,\ldots,v_n)) = (w_1,\ldots,w_n)$ where

$$\left[\sum_i u_i X_i(x), \sum_i v_i X_i(x)\right] = \sum_i w_i X_i(x).$$

As defined, $f(x) : \mathbb{R}^n \times \mathbb{R}^n \to \mathbb{R}^n$ is bilinear, i.e., $f(x) \in \mathrm{BL}(n,\mathbb{R})$. Using $\pi : G \to H$ from Corollary 3.6, we obtain a homomorphism $\Xi : G \to \mathrm{GL}(\mathrm{BL}(n,\mathbb{R}))$ defined so that $[\Xi(g)f](u,v) = \pi(g).f(\pi(g)^{-1}u, \pi(g)^{-1}v)$. Note, then, that $f(\rho(g)x) = [\Xi(g)f](x)$.

By ergodicity of the action, and arguments of Furstenberg [5] (see also [29, §3.2]), it is shown that $f$ is $G$-invariant, i.e., $\Xi(g)f = f$. Ergodicity of the action and continuity of $f$ thus force $f$ to be constant. Hence, $f(x) = f(\rho(g)x)$ for almost every $x \in M'$. Equivalently, there exist constants $C_{i,j}^k$ such that $[X_i, X_j] = \sum_k C_{i,j}^k X_k$, thereby producing the structure of a finite dimensional Lie algebra.

The proof is now completed by Proposition 3.3 in [4]. □

*Proof of Theorem* 1.4. The proof follows exactly as for Theorem 1.2; noting that Corollaries 3.5 and 3.6 hold in this situation as well. □

3.2. *Cartan actions of lattices.* In this subsection we consider actions of a lattice $\Gamma \subset G$ on a closed manifold $M$. In particular, we consider Cartan actions of such lattices [10].

*Definition* 3.7. Let $\mathcal{A}$ be a free abelian group and let $\{\gamma_1, \ldots, \gamma_n\}$ be a set of generators. A $C^r$ action $\rho : \mathcal{A} \times M \to M$ is called an *abelian Cartan action* if

1. for each $i$, $\rho(\gamma_i)$ is an Anosov diffeomorphism,

2. for each $i$, $\rho(\gamma_i)$ has a one dimensional strongest stable foliation $\mathcal{F}_i^{ss}$,

3. the tangential distributions $T\mathcal{F}_i^{ss}$ are pairwise transverse with internal direct sum $T\mathcal{F}_1^{ss} \oplus \cdots \oplus T\mathcal{F}_n^{ss} \cong TM$.

*Definition* 3.8. Let $\rho : \Gamma \times M \to M$ be an Anosov $C^r$ action. Then $\rho$ is called a *Cartan* (*lattice*) *action* if there is a subset of commuting hyperbolic elements $\{\gamma_1, \ldots, \gamma_n\} \subset \Gamma$ which generate an abelian subgroup $\mathcal{A}$ such that the restriction of $\rho|_\mathcal{A}$ is an abelian Cartan $C^r$ action on $M$.

*Remark* 3.9. Since the collection of strongest stable foliations is one dimensional for a Cartan action, it follows that a Cartan action must be multiplicity free.

*Example* 3.10 (Cartan actions). Consider the standard action of $\mathrm{SL}(n, \mathbb{Z})$ on the $n$-torus $\mathbb{T}^n \cong \mathbb{R}^n/\mathbb{Z}^n$. The super-rigidity representation $\pi$ obtained from Theorem 3.11 must be the standard representation of $\mathrm{SL}(n, \mathbb{R})$ on $\mathbb{R}^n$. Hence, each Lyapunov distribution is one dimensional, and, in fact, can be made the strongest stable distribution for some element in a fixed Cartan subgroup in $\mathrm{SL}(n, \mathbb{Z})$. It follows that the action is Cartan. See [6] for related results concerning actions of $\mathrm{SL}(n, \mathbb{Z})$ on a compact $n$ dimensional manifold. For additional examples of Cartan actions, and in particular on nilmanifolds, we refer the reader to [26].

Let $P \to M$ be the bundle of full frames over $M$, so that $P$ is a principal $\mathrm{GL}(n, \mathbb{R})$ bundle where $\mathrm{GL}(n, \mathbb{R})$ acts on the right on $P$. If $\Gamma$ acts on $M$ then there exists a natural lift via the derivative to a $\Gamma$ action on $P$. In [25], N. Qian established the following:

THEOREM 3.11. *Let $G$ be a connected semisimple Lie group with finite center and without compact factors such that each simple component of $G$ has $\mathbb{R}$-rank at least two. Let $\Gamma \subset G$ be a lattice. Let $M$ be a compact n dimensional smooth manifold with a measure $\mu$ taking positive values on open sets. Let $\rho$ be a $C^1$ Cartan action of $\Gamma$ on $M$ which is ergodic with respect to $\mu$. Then there exist*

1. *a subgroup $\Gamma_0 \subset \Gamma$ of finite index,*
2. *a $C^0$ section $\phi$ of the frame bundle $P \to M$, and*
3. *a homomorphism $\pi : \Gamma_0 \to \mathrm{GL}(n, \mathbb{R})$*

*such that with respect to the induced action of $\Gamma_0$ on $P \to M$, $\phi$ is a totally $\pi$-simple section, i.e., $\phi(\gamma x) = \gamma.\phi(x).\pi(\gamma)^{-1}$ for every $x \in M$, and every $\gamma \in \Gamma_0$.*

In fact, Qian proves this result holds under the weaker assumption that $\rho$ is weakly Cartan. We remark that in conjunction with Remark (d) following the proof of Theorem 1.1 in [25], this framing is actually Hölder. Further, because of the assumption that there exist Anosov elements in $\Gamma$, it follows that the measurable algebraic hull of the $\Gamma$ action on $P$ cannot contain compact factors, for otherwise the Anosov elements would admit zero Lyapunov exponents. Hence, the homomorphism $\pi$ extends to a homomorphism $\pi : G \to H$.

To apply our results from the previous section, we induce the $\Gamma$ action on $M$ to an action of $G$ on $X = \dfrac{G \times M}{\Gamma}$ (for ease of notation and without loss of generality, we may assume $\Gamma = \Gamma_0$). By construction, $X \to G/\Gamma$ is a fiber bundle with compact fibers $M$. Define actions of $G$ and $\Gamma$ on $G \times M$ so that

$$g(h, m) = (gh, m) \text{ and } \gamma(h, m) = (h\gamma^{-1}, \gamma m),$$

and $G$ and $\Gamma$ actions on $G \times P$ so that

$$g(h, p) = (gh, p) \text{ and } \gamma(h, p) = (h\gamma^{-1}, \gamma p).$$

Next, define $\Phi : G \times M \to G \times P$ so that

$$\Phi(g, m) = (g, \phi(m)\pi(g)^{-1}).$$

Then, for any $\gamma \in \Gamma$,

$$\begin{aligned}\Phi(\gamma(g, m)) &= \Phi(g\gamma^{-1}, \gamma m) = (g\gamma^{-1}, \phi(\gamma m)\pi(\gamma)\pi(g)^{-1}) \\ &= (g\gamma^{-1}, \gamma\phi(m)\pi(g)^{-1}) = \gamma\Phi(g, m),\end{aligned}$$

so that $\Phi$ is $\Gamma$-equivariant, and therefore defines a function $\Psi : X = \dfrac{G \times M}{\Gamma} \to Y = \dfrac{G \times P}{\Gamma}$. As a principal bundle over $X$, $Y$ is the bundle of frames tangent to the fiber for the bundle $X \to G/\Gamma$. Thus $\Psi$ is a Hölder section of $Y \to X$ which produces a framing for the fiber over every $x \in X$. Most importantly, $\Psi$ is also totally $\pi$-simple for every $g \in G$:

$$\Psi(g[h,m]) = \Psi([gh,m]) = [gh, \phi(m)\pi(h)^{-1}\pi(g)^{-1}] = g.\Psi([h,m])\pi(g)^{-1}.$$

This establishes an analog of Corollary 3.4 for the induced action of $G$ on $X$. Let $A \subset \Gamma$ be the maximal abelian subgroup such that the action $\rho|_A$ is abelian Cartan. If $F_x \subset X$ is the fiber through $x \in X$, then $F_x \cong M$ and the Hölder section $\Psi$ produces a Hölder framing for the fiber which can then be used to construct a Hölder Riemannian metric $\|\cdot\|_x$ for the fiber. The decomposition $TM = \oplus E_i$ induces a corresponding decomposition of the tangent space of $F_x$, and since $\Psi$ is totally $\pi$-simple, we have precise knowledge of the rates of expansion and contraction by elements in $A$, i.e., $\|av\|_{ax} = e^{\chi_i(\log a)}\|v\|_x$ for every $v \in E_i$ and every $a \in A$. The fact that the original $\Gamma$ action on $M$ is Cartan implies there exists an $A$-invariant trellis tangent to $F_x$ for the $\Gamma$ action on $X$. In particular, the $G$ action on $X$ is trellised with respect to $A$. Consequently, if we restrict our attention to the direction of the fibers, then Assumptions (A1), (A4) and (A5) from Section 2 hold. We also note that the ergodicity of $G$ on the finite volume space $X$ ensures that Assumptions (A2) and (A3) hold.

If $X$ is not compact, which occurs whenever $\Gamma$ is not cocompact in $G$, we cannot directly apply our arguments from Section 2. However, as we have just noted, if we restrict our attention to the fibers of $X \to G/\Gamma$, then Assumptions (A1) through (A5) from Section 2 still hold. With the help of Remark 2.13, we can easily adapt the argument used in proving Theorem 2.8 to hold in this situation, and, consequently, obtain $C^\infty$ information about the fibers of $X \to G/\Gamma$.

We have nearly proved the following result.

THEOREM 3.12. *Suppose $G$ is a connected semisimple Lie group of higher rank without compact factors such that each simple factor of $G$ has $\mathbb{R}$-rank at least 2. Let $\Gamma \subset G$ be a lattice and suppose that $\Gamma$ acts on a closed manifold $M$ so that the action is Cartan and volume-preserving for some smooth measure taking positive values on open sets. Let $A$ be a maximal $\mathbb{R}$-split Cartan of $G$, with $\{\chi_i\}$ the set of weights of $\pi$ with respect to $A$. Then there exist*

1. *a $C^\infty$ Riemannian metric, $\|\cdot\|$, on $M$, and*

2. *a $C^\infty$ decomposition $TM = \bigoplus E_i$*

*such that for every $v \in E_i$ and $a \in A$*

$$\|av\| = e^{\chi_i(\log a)}\|v\|.$$

*In particular, there exists a $C^\infty$ totally $\pi$-simple framing of $M$.*

*Proof.* It remains only to consider the $\mathbb{R}$-rank two case. In this case, we use the fact that inside the stable manifold, the one dimensional strongest stable foliation is $C^\infty$ with uniformly $C^\infty$ leaves. Hence we can apply Journé's result, Lemma 2.25, directly to conclude that the $\mathcal{F}_i$ vary $C^\infty$ along any $\mathcal{F}_j$ unless they have weights which are negative multiples of one another. However, in this case, we can appeal to Proposition 2.23 to establish the required smoothness of the $\mathcal{F}_i$ along $\mathcal{F}_j$. □

*Proof of Theorem* 1.5. The proof is essentially the same as that for Theorem 1.2. We use Theorem 3.12 to provide a smooth linearizing framing for the action on a subgroup of finite index, and then, exactly as in the proof of Theorem 1.2, we apply the arguments in [25, §3]. That $H$ is nilpotent follows from Proposition 3.13. To conclude that the action is conjugate to the standard algebraic action on $\tilde{\pi}_1(\tilde{M})/\pi_1(\tilde{M})$, we refer to [11, Cor. 2] to conclude that the action has a fixed point, and point out that any affine action with a fixed point must be conjugate to the standard action. □

*Proof of Corollary* 1.6. Any action sufficiently close to $\phi$ must also be volume-preserving and Cartan [17], and therefore, just as $\phi$, be $C^\infty$ conjugate to the standard algebraic action on $\tilde{\pi}_1(M')/\pi_1(M')$. □

Although the following result is presumably well-known, we provide a proof for the sake of completeness.

PROPOSITION 3.13. *Suppose there exists $\Phi \in \mathrm{Aff}(H)$ such that $\Phi$ is an Anosov diffeomorphism. Then $H$ is nilpotent.*

*Proof.* From [4], we can write $\Phi$ as a composition of an automorphism of $H$ with left multiplication by an element of $H$, i.e., $\Phi = L_g \circ \phi$ for $\phi$ some automorphism of $H$. Let $\mathfrak{h}$ be the Lie algebra of $H$ and define $\psi$ to be the automorphism of $\mathfrak{h}$ induced by $\mathrm{Ad}(g) \circ \phi_*$. Let $\psi_s$ be the semisimple component of the Jordan decomposition of $\psi$, also an automorphism of $\mathfrak{h}$. Since, with respect to a right invariant metric on $H$, $(L_g)_*\phi_*(v)$ and $(R_{g^{-1}})_*(L_g)_*\phi_*(v)$ have the same norm for any $v$, $\Phi$ being Anosov implies that $\psi$ and hence $\psi_s$ cannot have eigenvalues of modulus one.

Let $\mathfrak{s}$ be the solvradical of $\mathfrak{h}$. Then $\psi_s$ descends to an automorphism $\tilde{\psi}_s$ of $\mathfrak{h}/\mathfrak{s}$. Since $\mathfrak{h}/\mathfrak{s}$ is a semisimple Lie algebra, some finite power of $\tilde{\psi}_s$ must be the adjoint for some element of $\mathfrak{h}/\mathfrak{s}$, i.e., we may assume $\tilde{\psi}_s = \mathrm{Ad}(h + \mathfrak{s})$. But $\mathfrak{h}/\mathfrak{s}$ being semisimple, $\mathrm{Ad}(h + \mathfrak{s})$ must contain eigenvalues of modulus one. Since, the eigenvalues of $\tilde{\psi}_s$ are a subset of those of $\psi_s$, this yields a contradiction unless $\mathfrak{h}/\mathfrak{s}$ is trivial. We therefore may assume that $H$ is solvable.

To see that $\mathfrak{h}$ is nilpotent, we intend to show that $\mathfrak{h}$ is equal to its nilradical. For this purpose, we may, without loss of generality, consider the complexification of $\mathfrak{h}$, which we will also denote by $\mathfrak{h}$. Let $\mathfrak{n}$ be the nilradical of $\mathfrak{h}$, so that $[\mathfrak{h}, \mathfrak{h}] \subset \mathfrak{n}$. If $\mathfrak{h} \neq \mathfrak{n}$, then we may pick some $X \in \mathfrak{h}$ such that $X \notin \mathfrak{n}$ and $X$ is an eigenvector for $\psi_s$ with eigenvalue $\lambda$. Note that $|\lambda| \neq 1$. To see that $\mathfrak{h}$ is nilpotent, it will suffice to show that $\mathbb{R}X + \mathfrak{n}$ is a nilpotent ideal.

We shall now describe two distinct filtrations for $\mathfrak{n}$. First, we have the descending central series. Let $\mathcal{C}^0\mathfrak{n} = \mathfrak{n}$ and define $\mathcal{C}^i\mathfrak{n} = [\mathfrak{n}, \mathcal{C}^{i-1}\mathfrak{n}]$. Since $\mathfrak{n}$ is nilpotent there exists $k$ such that $\mathcal{C}^k\mathfrak{n} = 0$. Therefore we have the filtration

$$\mathfrak{n} = \mathcal{C}^0\mathfrak{n} \supset \mathcal{C}^1\mathfrak{n} \supset \cdots \supset \mathcal{C}^k\mathfrak{n} = 0.$$

Note that $[X, \mathcal{C}^i\mathfrak{n}] \subset \mathcal{C}^i\mathfrak{n}$. For the second filtration, when $|\lambda| > 1$, order the eigenvalues $\{\lambda_1, \ldots, \lambda_r\}$ of $\psi_s$ on $\mathfrak{n}$ so that $|\lambda_i| \leq |\lambda_j|$ for every $i < j$, and if $|\lambda| < 1$, then order them so that $|\lambda_i| \geq |\lambda_j|$ for every $i < j$. If $V_i$ is the eigenspace of $\psi_s$ with eigenvalue $\lambda_i$, then by defining $W_i = \bigoplus_{j=i+1}^{r} V_i$, we have

$$\mathfrak{n} = W_0 \supset W_1 \supset \cdots \supset W_r = 0.$$

We claim that $\mathbb{R}X + \mathfrak{n}$ is nilpotent with nilpotency degree $kr$. To see this consider some element in $\mathcal{C}^l(\mathbb{R}X + \mathfrak{n})$ with $l > kr$:

$$Y = [a_l X + N_l, [\cdots [a_2 X + N_2, a_1 X + N_1] \cdots],$$

where $a_i \in \mathbb{R}$ and $N_i \in \mathfrak{n}$. By expanding this expression, we can write $Y$ as a linear combination of terms of the form

$$y = \operatorname{ad}(Y_l)\operatorname{ad}(Y_{l-1})\cdots\operatorname{ad}(Y_2)(Y_1)$$

where either $Y_i = X$ or $Y_i \in \mathfrak{n}$ for every $i$. Since $[X, \mathcal{C}^i\mathfrak{n}] \subset \mathcal{C}^i\mathfrak{n}$, if at least $k$ of the $Y_i$'s lie in $\mathfrak{n}$, then $y = 0$. However, by construction of the $W_i$'s, we have $[X, W_i] \subset W_{i+1}$. Hence $(\operatorname{ad}(X))^r(N) = 0$ for any $N \in \mathfrak{n}$. So if there exists a string of $r$ consecutive $Y_i$'s all equal to $X$, then again $y = 0$. Since $l > kr$, one of these situations always occurs; i.e., we must have either at least $k$ of the $Y_i$'s belong to $\mathfrak{n}$ or that there exists some $j$ such that $Y_i = X$ for $i = j, j+1, \ldots, j+r-1$. In either case, we have that $y = 0$, forcing $Y = 0$.

We conclude that $\mathbb{R}X + \mathfrak{n}$ is nilpotent. Since it is clearly an ideal, we may conclude $\mathfrak{h} = \mathfrak{n}$. □

208 Monte Carlo Way, Danville, CA

University of Michigan, Ann Arbor, MI
*E-mail address*: spatzier@math.lsa.umich.edu